\documentclass[12pt,reqno]{amsart}

\usepackage[utf8]{inputenc}
\usepackage[T1]{fontenc}
\usepackage{amssymb}
\usepackage{mathtools}
\usepackage[dvipsnames]{xcolor}
\usepackage{bm}
\usepackage[shortlabels]{enumitem}
\usepackage{bbm}
\usepackage{hyperref}
\newtheorem{thm}{Theorem}[section]

\newtheorem{corollary}[thm]{Corollary}
\newtheorem{lemma}[thm]{Lemma}
\newtheorem*{lemma*}{Lemma}
\newtheorem*{remark*}{Remark}
\newtheorem{proposition}[thm]{Proposition}

\theoremstyle{definition}
\newtheorem{definition}[thm]{Definition}
\theoremstyle{remark}

\newtheorem{remark}[thm]{Remark}
\newtheorem{question}[thm]{Question}
\newtheorem{example}[thm]{Example}
\usepackage[left=3.5cm,right=3.5cm,top=3cm,bottom=3cm]{geometry}
\numberwithin{equation}{section}
\allowdisplaybreaks
\newcommand{\cI}{\ensuremath{\mathcal{I}}}
\newcommand{\cD}{\ensuremath{\mathcal{D}}}

\newcommand{\FF}{\ensuremath{\mathbb{F}}}
\newcommand{\RR}{\ensuremath{\mathbb{R}}}
\newcommand{\NN}{\ensuremath{\mathbb{N}}}

\newcommand{\CC}{\ensuremath{\mathbb{C}}}
\newcommand{\SSS}{\ensuremath{\mathbb{S}}}

\newcommand{\ee}{\ensuremath{\bm{e}}}
\newcommand{\yy}{\ensuremath{\bm{y}}}

\newcommand{\XX}{\ensuremath{\mathbb{X}}}

\newcommand{\ZB}{\ensuremath{\mathcal{Z}}}
\newcommand{\YY}{\ensuremath{\mathbb{Y}}}

\newcommand{\YB}{\ensuremath{\mathcal{Y}}}
\newcommand{\EB}{\ensuremath{\mathcal{E}}}
\newcommand{\Ind}{\ensuremath{\mathbbm{1}}}
\newcommand{\Id}{\ensuremath{\mathrm{Id}}}

\newcommand{\BB}{\ensuremath{\mathcal{B}}}

\newcommand{\kk}{\ensuremath{\bm{k}}}
\newcommand{\ku}{\ensuremath{\widehat{\bm{k}}}}

\newcommand{\cG}
{\ensuremath{\mathcal{G}}}

\DeclareMathOperator{\sgn}{sign}
\DeclareMathOperator{\spn}{span}
\DeclareMathOperator{\supp}{supp}

\newcommand{\kl}{\ensuremath{\tilde{\bm{k}}}}

\author[M. Berasategui]{Miguel Berasategui}
\address{Miguel Berasategui\\
	 UBA - Pab I, Facultad de Ciencias Exactas y Naturales\\
	Universidad de Buenos Aires\\
	(1428), Buenos Aires, Argentina}
\email{mberasategui@dm.uba.ar}

\author[P. M. Bern\'a]{Pablo M. Bern\'a}
\address{Pablo M. Bern\'a\\
	Universidad Tecnológica Atlántico Mediterráneo -- UTAMED\\
	Facultad de Empresa Digital, Tecnología y Derecho\\
	Málaga, España}
\email{pablomanuel.berna@utamed.es}

		\author[S. J. Dilworth]{Stephen J. Dilworth}
	\address{Stephen J. Dilworth, 	Department of Mathematics, University of South Carolina, Columbia, SC 29208, USA} 
	\email{dilworth@math.sc.edu}

\author[D. Kutzarova]{Denka Kutzarova}
\address{Denka Kutzarova, Department of Mathematics, University of Illinois, Urbana IL 61801, USA} 
\email{denka@illinois.edu}

\author[T. Oikhberg]{Timur Oikhberg}
\address{Timur Oikhberg, Department of Mathematics, University of Illinois, Urbana IL 61801, USA} 
\email{oikhberg@illinois.edu}

\begin{document}

\title[Characterizations of the unit vector basis of $\mathtt{c}_0$ or $\ell_p$]{Characterizations of the unit vector basis of $\mathtt{c}_0$ or $\ell_p$ among
minimal systems}

\thanks{M. Berasategui was supported by the Grants ANPCyT PICT 2018-04104 and CONICET PIP 11220200101609CO. Pablo M. Berná was supported by the Grant PID2022-142202NB-I00 (Agencia Estatal de Investigación, Spain). S. J. Dilworth was supported by the Simons Foundation Collaboration Grant No. 849142.}



\begin{abstract}
Given a minimal system $\EB$ for a Banach space $\XX$, Casazza, Dilworth,
Kutzarova and Motakis introduced a multiplication $\otimes$ on
$\spn\left(\EB\right)$ and studied the Schauder bases for which
$\left\Vert x\otimes y\right\Vert\approx \left\Vert x\right\Vert\left\Vert
y\right\Vert$, asking whether every such basis is equivalent to the canonical
unit vector basis of $\mathtt{c}_0$ or $\ell_p$. We continue this study in the
more general setting of minimal systems that are not assumed to be Schauder
bases, and we obtain positive answers under hypotheses that are formally weaker
than the ones considered so far: partial democracy instead of democracy or
bidemocracy, and truncation quasi-greediness instead of quasi-greediness. Along
the way we show that the above multiplicative estimate forces invariance under
translations, we characterize when such a minimal system is a Schauder basis, and we
give examples showing that the one-sided variants of the estimate behave very
differently.
\end{abstract}

\maketitle
%
%

\section{Introduction}\label{sectionintro}

Identifying the canonical unit vector basis of $\mathtt{c}_0$ and of $\ell_p$
$\left(1\le p<\infty\right)$ by means of intrinsic, isomorphic properties is a
classical theme in Banach space theory. The best known result in this direction
is due to Altshuler \cite{A1976}, who proved that if $\EB$ is a symmetric basis
for $\XX$ and both $\XX$ and $\overline{\spn\left(\EB^*\right)}$ have a unique
symmetric basic sequence up to equivalence, then $\EB$ is equivalent to the
canonical unit vector basis of $\mathtt{c}_0$ or of $\ell_p$ for some
$1\le p<\infty$. Altshuler's theorem has recently been extended to subsymmetric
bases by Casazza, Dilworth, Kutzarova and Motakis \cite{CDKM2022}.

In the theory of greedy approximation with respect to bases, symmetry is
usually far too strong a requirement, and it is customary to work instead with
the much weaker notions that arise from the study of the Thresholding Greedy
Algorithm of Konyagin and Temlyakov \cite{KT1999}: quasi-greediness,
democracy, bidemocracy, almost greediness, and so on (see
\cite{AK2016,T2011,W2000} and the references therein). It is therefore natural
to ask for versions of Altshuler's theorem in which symmetry is replaced by
hypotheses of greedy type. This is precisely the programme initiated in
\cite{CDKM2023}, and the present paper is a contribution to it.

\subsection{The multiplication \texorpdfstring{$\otimes$}{tensor}}
Let $\EB$ be a minimal system for a Banach space $\XX$ and, for each
$x, y\in \spn\left(\EB\right)$, define 
\begin{align*}
x\otimes y=x\otimes y \left[\EB,\XX\right]:=\sum_{n=1}^{\max(\lbrace{0\rbrace}\cup \supp(x))}\ee_n^*\left(x\right)\sum_{m=1}^{\max(\lbrace{0\rbrace}\cup \supp(y))}\ee_m^*\left(y\right)\ee_{\left(n-1\right)m\left(y\right)+m},
\end{align*}
and let $x^2:=x\otimes x$ (note that $x\otimes 0=0\otimes x=0$). Thus $\left(\spn\left(\EB\right),\otimes\right)$ is a
non-commutative semigroup with unit $\ee_1$; the multiplication is, however,
not continuous.

If $\EB$ is equivalent to the canonical unit vector basis of $\mathtt{c}_0$ or
of $\ell_p$, then there is $K\ge 1$ such that
\begin{align}
&\frac{1}{K}\left\Vert x\right\Vert \left\Vert y\right\Vert\le \left\Vert x\otimes y\right\Vert\le K \left\Vert x \right\Vert \left\Vert y\right \Vert \label{p**}
\end{align}
for all $x, y \in \spn\left(\EB\right)$. In \cite{CDKM2023}, the authors
introduced and studied the Schauder bases satisfying \eqref{p**}, and posed the
following question, which is the main motivation for this note.

\begin{question}\label{questionCDKM}\cite[Question 1]{CDKM2023}
Let $\EB$ be a seminormalized Schauder basis for $\XX$ satisfying \eqref{p**}.
Is $\EB$ equivalent to the canonical unit vector basis of $\mathtt{c}_0$ or of
$\ell_p$ for some $1\le p<\infty$?
\end{question}

Among their main results, they proved that Question~\ref{questionCDKM} has a
positive answer whenever the basis is bidemocratic, almost greedy, or invariant
under spreading (see \cite[Theorems 3.1, 5.1 and 6.1]{CDKM2023}).

\subsection{Our contribution}
In this paper we continue the study of the above mentioned property,
strengthening some of the known results and extending some of them to arbitrary
minimal systems. Our work departs from \cite{CDKM2023} in three directions.

\emph{(a) From Schauder bases to minimal systems.} We do not assume a priori
that $\EB$ is a Schauder basis, nor even that it is a Markushevich basis. This
is not a cosmetic generalization: on the one hand, minimal systems which are
not Schauder bases occur naturally in greedy approximation theory; on the
other hand, Example~\ref{examplenotproduct} exhibits a normalized,
$1$-translation invariant minimal system that satisfies the upper estimate
\eqref{p*} and the lower estimate for powers \eqref{p*l-+} with constant
$K=1$ and is \emph{not} a Markushevich basis. What we prove instead is that the
Schauder property can be \emph{recovered} rather than assumed: by
Lemma~\ref{lemmaSchauder3}, a minimal system satisfying \eqref{p**} is a
Schauder basis if and only if $\EB^*$ is bounded, if and only if the
auxiliary norm $N_{\cD}$ of Definition~\ref{definitiondifferences} is not
equivalent to the norm of the space, and the failure of these conditions forces
an extremely fast growth of the biorthogonal functionals together with
$\varphi\left(n\right)\gtrsim n$.

\emph{(b) Weaker forms of the multiplicative estimate.} Following
\cite{CDKM2023}, we split \eqref{p**} into the one-sided and ``power''
variants \eqref{p*}--\eqref{p*l-} listed in
Section~\ref{subsectionproperties}, and we keep track, in each statement, of
which of them is actually needed. Remark~\ref{remarkmax} and
Example~\ref{examplenotproduct} show that the resulting hierarchy is genuinely
finer than \eqref{p**}.

\emph{(c) Weaker forms of the greedy-type hypotheses.} Democracy is replaced by
\emph{partial democracy}, and quasi-greediness by \emph{truncation
quasi-greediness}, both of which make sense for arbitrary minimal systems and
are formally weaker than their classical counterparts. The reason why this
costs us nothing is a structural phenomenon that seems to be of independent
interest: \eqref{p**} \emph{forces invariance under translations}
(Lemma~\ref{lemmatranslates}). Under translation invariance, democracy,
conservativeness and partial democracy turn out to be equivalent
(Lemma~\ref{lemmademoc0}), and a weak bidemocracy estimate for intervals holds
automatically (Lemma~\ref{lemmabidemlite}); combining the two we obtain that a
translation invariant Schauder basis for which both $\EB$ and $\EB^*$ are
partially democratic is bidemocratic (Corollary~\ref{corollarybidemtranslate}).
In a similar vein, Lemma~\ref{lemmauncondparam2} shows that a minimal system
satisfying only \eqref{p*-} and \eqref{p*l-} is either unconditional with
constant $K^2$ or so conditional that its conditionality parameters $\kk_m$ grow faster than a power $m^{\alpha}$, for some $\alpha>0$; since truncation quasi-greedy seminormalized systems
have conditionality parameters of order $\log m$, the second alternative is
ruled out and unconditionality follows
(Corollaries~\ref{corollaryuncond12} and~\ref{corollaryuncond2}).

Our main result, Proposition~\ref{propositionnewcharac}, collects six sets of
hypotheses --- each of them formally weaker than bidemocracy or almost
greediness --- under which a minimal system satisfying \eqref{p**} must be a
Schauder basis equivalent to the canonical unit vector basis of $\mathtt{c}_0$
or of $\ell_p$ for some $1\le p<\infty$. We complement it with
Proposition~\ref{propositionsubseqtoellp}, where invariance under spreading is
replaced by the requirement that every subsequence of $\EB$ be translation
invariant; the equivalence of these two conditions, established in
Lemma~\ref{lemmatranslationtosubsym}, may also be of some independent interest.

\subsection{Organization of the paper}
Section~\ref{sectionnotation} fixes the notation and the terminology, and
records the elementary properties of the multiplication $\otimes$ that will be
used throughout. Section~\ref{sectionbasic} contains some basic consequences of
\eqref{p**} and of its weaker variants. Section~\ref{sectiontranslations} is
devoted to translation invariance and to its relation with invariance under
spreading. In Section~\ref{sectionschauder} we determine when a minimal system
satisfying \eqref{p**} is a Schauder basis. Section~\ref{sectiontqg} deals with
truncation quasi-greediness and unconditionality, and
Section~\ref{mainresults} contains the proof of our main result and the above mentioned complement.

\section{General notation and terminology}\label{sectionnotation}

Throughout the paper we use standard Banach space notation and terminology as
in \cite{AK2016,LT1977}, and we refer the reader to \cite{T2011} for the
background on greedy approximation. This section collects, for the reader's
convenience, the conventions that will be used from now on and the
definitions of the properties involved in our results. All of the notation
below extends from Banach spaces to normed spaces, and to general vector
spaces, whenever this makes sense.

\subsection{Minimal systems}\label{subsectionminimal}
Unless otherwise specified, $\XX$ denotes a Banach space over the field $\FF$,
where $\FF=\RR$ or $\FF=\CC$, and $\EB=\left(\ee_n\right)_{n\in \NN}$ denotes a
\emph{minimal system} for $\XX$; that is, a sequence whose linear span is dense
in $\XX$ and for which there is a sequence
$\EB^*=\left(\ee_n^*\right)_{n\in \NN}\subset \XX^*$ with
$\ee_n^*\left(\ee_m\right)=\delta_{n,m}$ for all $m,n\in \NN$. We call $\EB^*$
the \emph{dual system} of $\EB$, and we call \emph{bidual system} of $\EB$ the
sequence $\EB^{**}=\left(\ee_n^{**}\right)_{n\in \NN}$ of functionals in
$\overline{\spn\left(\EB^*\right)}^{*}$ given by duality.

The system $\EB$ is \emph{normalized} if $\left\Vert \ee_n\right\Vert=1$ for
every $n\in \NN$, and \emph{seminormalized} if
$0<\inf_{n\in \NN}\left\Vert \ee_n\right\Vert\le \sup_{n\in\NN}\left\Vert \ee_n\right\Vert<\infty$;
we write $M\left[\EB\right]:=\sup_{n\in\NN}\left\Vert \ee_n\right\Vert$, where
the value $\infty$ is allowed. We say that $\EB$ is a \emph{Markushevich basis}
(M-basis, for short) if $\EB^*$ is total, i.e., if $\ee_n^*\left(x\right)=0$
for every $n\in \NN$ entails $x=0$, and that it is a \emph{norming} M-basis if,
moreover, $\spn\left(\EB^*\right)$ is a norming subspace of $\XX^*$, i.e., if
there is $C\ge 1$ such that
\begin{align*}
&\left\Vert x\right\Vert\le C\sup\left\lbrace \left\vert x^*\left(x\right)\right\vert: x^*\in \SSS_{\XX^*}\cap \spn\left(\EB^*\right)\right\rbrace&&\forall x\in \XX.
\end{align*}
Schauder bases and unconditional bases are discussed separately in
Subsection~\ref{subsectionschauderbases}.

Every piece of notation we introduce depends on the pair $\left(\EB,\XX\right)$,
and we display this dependence in brackets whenever there is any risk of
confusion, writing for instance $\supp\left[\EB,\XX\right]\left(x\right)$ or
$\varphi\left[\EB,\XX\right]\left(m\right)$; otherwise we simply write
$\supp\left(x\right)$ or $\varphi\left(m\right)$. The corresponding objects
associated with the pair $\left(\EB^*,\XX^*\right)$ are denoted with an
asterisk, as in $\varphi^*\left(m\right)=\varphi\left[\EB^*,\XX^*\right]\left(m\right)$.

Finally, $\mathtt{c}_{00}$ stands for the vector space of all finitely supported
scalar sequences, $\mathtt{c}_0$ and $\ell_p$ for the usual sequence spaces, and
the \emph{canonical vector system} of a space of sequences is its sequence of
coordinate vectors.

\subsection{Sets, sequences, and domination}\label{subsectionsets}
Given a set $A$ and functions $f,g: A\rightarrow \RR$, we say that $g$
\emph{dominates} $f$, and we write $f\lesssim g$, if there is $C>0$ such that
$f\left(a\right)\le Cg\left(a\right)$ for every $a\in A$; if we wish to keep
track of the constant, we say that $g$ \emph{$C$-dominates} $f$. If
$f\lesssim g$ and $g\lesssim f$ we write $f\approx g$ and say that $f$ and $g$
are \emph{equivalent}, or \emph{$C$-equivalent} if we wish to record the
constants. For functions with values in a normed space, domination and
equivalence refer to the corresponding norms.

We write $\left\vert A\right\vert$ for the cardinality of a set $A$,
$A^{<\infty}$ for the set of its finite subsets and, given $m\in \NN_0$,
$A^{\left(m\right)}$ and $A^{\le m}$ for the sets of subsets of $A$ of
cardinality equal to $m$ and at most $m$, respectively; we use similar notation
for $<$, $>$ and $\ge$. If $A, B\subset \RR$ and $t\in \RR$, then $A<B$ means
that $a<b$ for all $a\in A$ and all $b\in B$, whereas $t<A$ (resp. $A<t$) means
$\left\lbrace t\right\rbrace<A$ (resp. $A<\left\lbrace t\right\rbrace$); note
that $\emptyset<A$ and $A<\emptyset$ for every $A\subset \RR$. We also set
$A_{\ge t}:=A\cap\left[t,+\infty\right)$ (and similarly for $\le$, $<$ and
$>$) as well as $a+A:=\left\lbrace a+x: x\in A\right\rbrace$ and
$aA:=\left\lbrace ax: x\in A\right\rbrace$ for $a\in \RR$. The unit sphere of a
normed space $A$ is denoted by $\SSS_A$. By $\cI$ we denote the set of finite
intervals of $\NN$; given $n,m\in \NN$, we write
$I_{n,m}:=\left\lbrace n,\dots,m\right\rbrace$, with the convention that
$I_{n,m}=\emptyset$ if $m<n$, and we abbreviate $I_m:=I_{1,m}$.

Given $\alpha=\left(\alpha_n\right)_{n\in \NN}\in \FF^{\NN}$, its
\emph{primitive sequence} is
$s\left(\alpha\right)=\left(s_n\left(\alpha\right)\right)_{n\in \NN}$, where
$s_n\left(\alpha\right):=\sum_{k=1}^{n}\alpha_k$ for $n\in \NN$. For $w\in \FF$
we set $\sgn\left(w\right):=w/\left\vert w\right\vert$ if $w\not=0$, and
$\sgn\left(w\right):=1$ if $w=0$.

\subsection{Supports, indicator vectors and projections}\label{subsectionsupports}
Let $\EB$ be a minimal system for $\XX$. The \emph{support} of $x\in \XX$ with
respect to $\EB$ is the set
$\supp\left(x\right):=\left\lbrace n\in \NN: \ee_n^*\left(x\right)\not=0\right\rbrace$,
and we set
\begin{align*}
&m\left(x\right):=\sup\left(\left\lbrace 0\right\rbrace \cup  \supp\left(x\right)\right),
&&l\left(x\right):=\min\left(\left\lbrace m\left(x\right)\right\rbrace  \cup  \supp\left(x\right)\right),
\end{align*}
where the value $\infty$ is allowed for $m\left(x\right)$. We also put
$\varepsilon\left(x\right):=\left(\sgn\left(\ee_n^*\left(x\right)\right)\right)_{n\in \NN}$
and, for $1\le p\le \infty$,
$\left\Vert x\right\Vert_{\ell_p}:=\left(\sum_{n\in \NN}\left\vert \ee_n^*\left(x\right)\right\vert^p\right)^{1/p}\in\left[0,+\infty\right]$.

Given $A\in \NN^{<\infty}$, $A\subset B\subset \NN$ and
$\varepsilon=\left(\varepsilon_n\right)_{n\in B}\in \SSS_{\FF}^{B}$, the
associated \emph{indicator vectors} and \emph{indicator functionals} are
\begin{align*}
&\Ind_{\varepsilon, A}:=\sum_{n\in A}\varepsilon_n\ee_n,&&\Ind_{A}:=\sum_{n\in A}\ee_n,&&&\Ind_{\varepsilon, A}^*:=\sum_{n\in A}\varepsilon_n\ee_n^*,&&&&\Ind_{A}^*:=\sum_{n\in A}\ee_n^*.
\end{align*}
Given $x\in \XX$, $B\in \NN^{<\infty}$, $B\subset D\subset \NN$,
$\varepsilon\in \SSS_{\FF}^{D}$ and $\alpha=\left(\alpha_n\right)_{n\in\NN}\in \FF^{\NN}$,
we consider the coordinate projections and the diagonal multipliers
\begin{align*}
&P_B\left(x\right):=\sum_{n\in B}\ee_n^*\left(x\right)\ee_n, &&P_{B,\varepsilon}\left(x\right):=\sum_{n\in B}\varepsilon_n\ee_n^*\left(x\right)\ee_n, &&&\alpha \cdot x:=\sum_{n\in \NN}\alpha_n\ee_n^*\left(x\right)\ee_n,
\end{align*}
the last one being defined for $x\in \spn\left(\EB\right)$. The operators
$P_B$ and $P_{B,\varepsilon}$ are bounded and finitely supported, so no
convergence issue arises; we call $\left(P_B\right)_{B\in\NN^{<\infty}}$ the
\emph{coordinate projections} associated with $\EB$. The \emph{partial sum
operators} and their complements are
\begin{align*}
&S_{n}\left(x\right):=P_{I_n}\left(x\right)&&\text{and}&&&S_{n,c}\left(x\right):=x-S_n\left(x\right)&&&&\left(n\in\NN\right),
\end{align*}
with the convention $S_0\left(x\right):=0$ and $S_{0,c}\left(x\right):=x$.
Finally, for $n\in \NN_0$ we put $\EB_{>n}:=\left(\ee_k\right)_{k>n}$ and
$\XX_{>n}:=\overline{\spn\left(\EB_{>n}\right)}$, and we use similar notation
for $<$, $\le$ and $\ge$.

Given $0<t\le 1$ and $x\in \XX$, a set $A\subset \NN$ is a \emph{$t$-greedy set}
of $x$ with respect to $\EB$ if
\begin{align*}
&\left\vert \ee_n^*\left(x\right)\right\vert \ge t \left\vert \ee_k^*\left(x\right)\right\vert &&\forall n\in A, \forall k\not\in A,
\end{align*}
and it is a \emph{greedy set} of $x$ if it is a $1$-greedy set. For each
$m\in \NN_0$, $\cG\left(x,m,t\right)$ denotes the set of $t$-greedy sets of $x$
of cardinality $m$, whereas $\cG\left(x,t\right)$ denotes the set of all finite
$t$-greedy sets of $x$. A map $\rho: \NN\rightarrow \NN$ is a \emph{greedy
ordering} for $x\in \XX$ if
$\left\lbrace \rho\left(j\right): 1\le j\le m\right\rbrace\in \cG\left(x,m,1\right)$
for every $m\in \NN$.

\subsection{Schauder bases and unconditionality}\label{subsectionschauderbases}
A minimal system $\EB$ for $\XX$ is a \emph{Schauder basis} if every $x\in \XX$
is the sum of a norm-convergent series of the form $\sum_{n\in\NN}a_n\ee_n$
with $\left(a_n\right)_{n\in \NN}\in \FF^{\NN}$. In that case the coefficients
are uniquely determined by $a_n=\ee_n^*\left(x\right)$, so that
\begin{align}
&x=\sum_{n\in \NN}\ee_n^*\left(x\right)\ee_n=\lim_{n\to \infty}S_n\left(x\right)&&\forall x\in \XX, \label{schauderexpansion}
\end{align}
and $\EB$ is in particular an M-basis. Since $\spn\left(\EB\right)$ is dense in
$\XX$ and $S_n\left(y\right)=y$ for every $y\in \spn\left(\EB\right)$ and every
$n\ge m\left(y\right)$, a routine argument shows that $\EB$ is a Schauder basis
if and only if the partial sum operators are uniformly bounded. When this is
the case, the \emph{basis constant} of $\EB$ and its complementary counterpart
are
\begin{align*}
&K_{b}=K_{b}\left[\EB,\XX\right]:=\sup_{n\in \NN_0}\left\Vert S_{n}\right\Vert
&&\text{and}
&&&K_{b,c}=K_{b,c}\left[\EB,\XX\right]:=\sup_{n\in \NN_0}\left\Vert S_{n,c}\right\Vert,
\end{align*}
and they satisfy $1\le K_b$, $K_{b,c}\le 1+K_b$ and $K_b\le 1+K_{b,c}$ because
$S_{n,c}=\Id_{\XX}-S_n$. Following the usual terminology, $\EB$ is
\emph{monotone} if $K_b=1$, and \emph{bimonotone} if $K_b=K_{b,c}=1$; every
Schauder basis admits an equivalent renorming for which it is bimonotone.
Moreover, \eqref{schauderexpansion} implies that $\spn\left(\EB^*\right)$ is
$K_b$-norming, so every Schauder basis is a norming M-basis, and
$1\le \left\Vert \ee_n\right\Vert\left\Vert \ee_n^*\right\Vert\le 2K_b$ for
every $n\in \NN$; in particular, $\EB^*$ is seminormalized whenever $\EB$ is.

A Schauder basis $\EB$ is \emph{unconditional} if the series in
\eqref{schauderexpansion} converges unconditionally for every $x\in \XX$.
As is well known, this happens if and only if the coordinate projections
$\left(P_A\right)_{A\in \NN^{<\infty}}$ are uniformly bounded, and also if and
only if the diagonal multipliers $x\mapsto \alpha\cdot x$ are uniformly bounded
for $\alpha\in \SSS_{\FF}^{\NN}$. The two associated constants are the
\emph{suppression unconditional constant} and the \emph{unconditional
constant},
\begin{align*}
&K_{su}=K_{su}\left[\EB,\XX\right]:=\sup\left\lbrace \left\Vert P_A\right\Vert: A\in \NN^{<\infty}\right\rbrace,\\
&K_{u}=K_{u}\left[\EB,\XX\right]:=\sup\left\lbrace \left\Vert \alpha\cdot x\right\Vert: x\in \SSS_{\XX}\cap \spn\left(\EB\right), \alpha\in \SSS_{\FF}^{\NN}\right\rbrace,
\end{align*}
and we say that $\EB$ is \emph{$C$-suppression unconditional} (resp.
\emph{$C$-unconditional}) if $K_{su}\le C$ (resp. $K_u\le C$). One always has
$K_b\le K_{su}\le K_u$ and, by convexity, $K_u\le 2K_{su}$ if $\FF=\RR$ and
$K_u\le 4K_{su}$ if $\FF=\CC$, so that the two constants are equivalent and
either of them may be used to quantify unconditionality. As with bimonotonicity,
every unconditional basis can be renormed so as to become $1$-unconditional.

Since our hypotheses are placed on minimal systems that are not assumed beforehand to be Schauder bases, we record the following elementary observation, which will be used without further reference.

\begin{remark}\label{remarkboundedprojections}
Let $\EB$ be a minimal system for $\XX$. If
$\sup\left\lbrace \left\Vert P_A\right\Vert: A\in \NN^{<\infty}\right\rbrace=:C<\infty$,
then $\EB$ is an unconditional Schauder basis with $K_{su}=C$ and $K_b\le C$.
\end{remark}

\subsection{Fundamental functions and conditionality parameters}\label{subsectionfundamental}
The \emph{fundamental function} $\varphi$, the \emph{upper democracy function}
$\varphi_u$ and the function $\lambda$ associated with $\EB$ are given, for
$m\in \NN$, by
\begin{align*}
&\varphi\left(m\right):=\sup_{A\in\NN^{\le m}}\left\Vert \Ind_{A}\right\Vert,
&&\varphi_u\left(m\right):=\sup_{\substack{A\in\NN^{\le m}\\\varepsilon\in \SSS_{\FF}^{A}}}\left\Vert \Ind_{\varepsilon, A}\right\Vert,
&&&\lambda\left(m\right):=\left\Vert \Ind_{I_m}\right\Vert,
\end{align*}
and we set $\lambda\left(0\right)=\varphi\left(0\right)=\varphi_u\left(0\right)=0$
for convenience. The \emph{conditionality parameters} of $\EB$ are, for
$m\in \NN$,
\begin{align*}
&\kk_m:=\sup \left\lbrace \left\Vert P_A\right\Vert: A\in \NN^{\le m}\right\rbrace;\\
&\ku_m:=\sup \left\lbrace  \left\Vert \varepsilon\cdot x\right\Vert: x\in \SSS_{\XX_{\le m}},  \varepsilon \in \SSS_{\FF}^{\NN} \right\rbrace;\\
&\kl_m:=\sup \left\lbrace \left\Vert P_A\right\Vert: A\subset I_m \right\rbrace.
\end{align*}
These are the finite-dimensional truncations of the constants of
Subsection~\ref{subsectionschauderbases}: one has
\begin{align*}
&\sup_{m\in \NN}\kk_m=\sup_{m\in \NN}\kl_m=K_{su}&&\text{and}&&&\sup_{m\in \NN}\ku_m=K_u,
\end{align*}
where the value $\infty$ is allowed. Consequently, and by
Remark~\ref{remarkboundedprojections}, a minimal system $\EB$ is an
unconditional Schauder basis if and only if one (equivalently, any)  of the
sequences $\left(\kk_m\right)_{m\in \NN}$, $\left(\ku_m\right)_{m\in \NN}$ and
$\left(\kl_m\right)_{m\in \NN}$ is bounded, and the rate at which these
sequences grow measures how far $\EB$ is from being unconditional.

\begin{remark}\label{remarkku}Note that by convexity, $\ku_m\le 2\kl_m$ if $\FF=\RR$, and $\ku_m\le 4\kl_m$ if $\FF=\CC$; this is the finite-dimensional counterpart of the estimate $K_u\le 2K_{su}$ (resp. $K_u\le 4K_{su}$) recorded in Subsection~\ref{subsectionschauderbases}.
\end{remark}

\subsection{Greedy-type properties}\label{subsectiongreedytype}
We will use the following notions from the literature on greedy approximation
(see \cite{AK2016,T2011}), all of which make sense for an arbitrary minimal
system. Below, $\EB$ is a minimal system for $\XX$ and $C>0$.

\begin{definition}\label{definitiongreedytype}
We say that $\EB$ is
\begin{enumerate}[\rm (i)]
\item\label{defsuppuncondcc} \emph{$C$-suppression unconditional for constant coefficients} if
\begin{align*}
&\left\Vert \Ind_{\varepsilon, A}\right\Vert\le C\left\Vert\Ind_{\varepsilon, B}\right\Vert&&\forall A\subset B\in \NN^{<\infty},\forall \varepsilon\in \SSS_{\FF}^{B};
\end{align*}
\item\label{defuncondcc} \emph{$C$-unconditional for constant coefficients} if
\begin{align*}
&\left\Vert \Ind_{\varepsilon, A}\right\Vert\le C\left\Vert\Ind_{\varepsilon', A}\right\Vert&&\forall A\in \NN^{<\infty}, \forall \varepsilon, \varepsilon'\in \SSS_{\FF}^{A};
\end{align*}
\item\label{defsuperdem} \emph{$C$-superdemocratic} if
\begin{align*}
&\left\Vert \Ind_{\varepsilon, A}\right\Vert\le C\left\Vert\Ind_{\varepsilon', B}\right\Vert&&\forall A, B\in \NN^{<\infty}\!: \left\vert A\right\vert\le \left\vert B\right\vert, \forall \varepsilon, \varepsilon'\in \SSS_{\FF}^{A\cup B},
\end{align*}
and \emph{$C$-democratic} if the above holds for $\varepsilon_n=\varepsilon'_n=1$ for all $n\in A\cup B$;
\item\label{defcons} \emph{$C$-conservative} if
\begin{align*}
&\left\Vert \Ind_{A}\right\Vert\le C\left\Vert \Ind_{ B}\right\Vert&&\forall A<B\in \NN^{<\infty}, \left\vert A\right\vert\le \left\vert B\right\vert;
\end{align*}
\item\label{defpardem} \emph{$C$-partially democratic} if for each $A\in\NN^{<\infty}$ there is $A\subset D\in\NN^{<\infty}$ such that
\begin{align*}
&\left\Vert \Ind_{A}\right\Vert\le C\left\Vert\Ind_{B}\right\Vert&&\forall B\in \NN^{<\infty}\!: B\cap D=\emptyset, \left\vert A\right\vert\le \left\vert B\right\vert;
\end{align*}
\item\label{defbidem} \emph{$C$-bidemocratic} if
\begin{align*}
&\varphi\left(m\right)\varphi^*\left(m\right)\le C m &&\forall m\in \NN;
\end{align*}
\item\label{defqg} \emph{$C$-quasi-greedy} if
\begin{align*}
&\left\Vert P_A\left(x\right)\right\Vert\le C\left\Vert x\right\Vert &&\forall x\in \XX, \forall A\in \cG\left(x,1\right);
\end{align*}
\item\label{deftqg} \emph{$C$-truncation quasi-greedy} if
\begin{align*}
&\min_{n\in A}\left\vert \ee_n^*\left(x\right)\right\vert \left\Vert\Ind_{\varepsilon\left(x\right), A}\right\Vert\le C\left\Vert x\right\Vert&&\forall x\in \XX, \forall A\in \cG\left(x,1\right).
\end{align*}
\end{enumerate}
As usual, we drop the constant from the terminology when we do not need to keep
track of it.
\end{definition}

Every quasi-greedy minimal system is truncation quasi-greedy, and the converse
does not hold; this is the reason why hypothesis \ref{deftqg} is weaker than
hypothesis \ref{defqg} in our statements.

We will also need the following regularity condition. A sequence of positive
numbers $\left(r_n\right)_{n\in \NN}$ has the \emph{upper regularity property}
(URP) if there is $b\in \NN_{\ge 2}$ such that $r_{bn}\le \frac{b}{2}r_n$ for
every $n\in \NN$, and a minimal system has the URP if its fundamental function
does.

\subsection{Translations and invariance under spreading}\label{subsectiontranslations}
The following operators, and the invariance properties associated with them,
play a central role throughout the paper.

\begin{definition}\label{definitiontranslate}Let $\EB$ be a minimal system for $\XX$. For each $k\in \NN_0$ and each $x\in \spn\left(\EB\right)$, we define the $k$-translation $T_k:\spn\left(\EB\right)\rightarrow \spn\left(\EB_{>k}\right)$ by
\begin{align*}
&T_k\left(x\right)=T_k\left[\EB,\XX\right]\left(x\right):=\sum_{n\in \NN}\ee_n^*\left(x\right)\ee_{n+k}.
\end{align*}
Given $K>0$, we say that $\EB$ is:
\begin{itemize}
\item $K$-left-translation invariant if
\begin{align*}
& \left\Vert x\right\Vert\le K \left\Vert T_k\left(x\right)\right\Vert&&\forall x\in \spn\left(\EB\right), k\in \NN_0.
\end{align*}
\item $K$-right-translation invariant if
\begin{align*}
&\left\Vert T_k\left(x\right)\right\Vert\le K \left\Vert x\right\Vert&&\forall x\in \spn\left(\EB\right), k\in \NN_0.
\end{align*}
\item $K$-translation invariant if it is $K$-left-translation invariant and $K$-right-translation invariant.
\end{itemize}
We say that $\EB$ is left-translation invariant (resp. right-) if there is $K>0$ for which $\EB$ is $K$-left(resp. right)-translation invariant.
\end{definition}

\begin{definition}\label{weaklysubsym}\cite{CDKM2023} A minimal system $\EB=\left(\ee_n\right)_{n\in \NN}$  for $\XX$ is \emph{invariant under spreading} with constant $C$ if every subsequence of $\EB$ is $C$-equivalent to $\EB$.
\end{definition}

\subsection{{Weaker variants of the multiplicative estimate}}\label{subsectionproperties}

Recall from Section~\ref{sectionintro} the definition of the
multiplication $\otimes$ on $\spn\left(\EB\right)$ and of the property
\eqref{p**}. Following \cite{CDKM2023}, we will also study weaker variants
thereof; in particular, we will consider the following properties:
\begin{align}
&\left\Vert x\otimes y\right\Vert\le K \left\Vert x \right\Vert \left\Vert y\right \Vert; \label{p*}\\
&\left\Vert x^n\right\Vert\le K^{n-1} \left\Vert x \right\Vert^n &&\forall n\in \NN;  \label{p*-+} \\
&\left\Vert x^2\right\Vert\le K \left\Vert x \right\Vert^2 ;  \label{p*-} \\
&\frac{1}{K}\left\Vert x\right\Vert \left\Vert y\right\Vert\le \left\Vert x\otimes y\right\Vert; \label{p*l}\\
&\frac{1}{K^{n-1}}\left\Vert x\right\Vert^n \le \left\Vert x^n\right\Vert &&\forall n\in \NN;\label{p*l-+}\\
&\frac{1}{K}\left\Vert x\right\Vert^2 \le \left\Vert x^2\right\Vert,\label{p*l-}
\end{align}
where $x$ and $y$ might be any elements in $\spn\left(\EB\right)$.  

The following are elementary properties of the function $\otimes$. Most of them are easy to check, so the proofs are left to the reader. 

\begin{remark}\label{remarksumexponents}For each $x\in \spn\left(\EB\right)$ and each $k, m\in \NN$, 
\begin{align*}
x^{k}\otimes x^{m}=& x^{k+m}. 
\end{align*}
\end{remark}

\begin{remark}\rm \label{remarktranslate} For each $x\in \spn\left(\EB\right)$ and each $k\in \NN$, 
\begin{align*}
\ee_k\otimes x=&T_{\left(k-1\right) m\left(x\right)}\left(x\right). 
\end{align*}
\end{remark}

\begin{remark}\label{remarkdistrib} \rm If  $m\left(x\right)=m\left(y\right)=m\left(x+y\right)$, then for all $z\in \spn\left(\EB\right)$, 
\begin{align*}
&z\otimes \left(x+y\right)=z\otimes x+z\otimes y&&\text{and}  &&&\left(x+y\right)\otimes z=x\otimes z+y\otimes z. 
\end{align*}
\end{remark}

\begin{remark}\label{remarkiterate}
If $x,y\in \spn\left(\EB\right)\setminus \left\lbrace 0\right\rbrace$ and $j\in \NN$, then 
$j\in \supp\left(x\otimes y\right)$ if and only if 
\begin{align*}
0\not=\ee_j^*\left(x\otimes y\right)=\sum_{n=1}^{m\left(x\right)}\sum_{m=1}^{m\left(y\right)}\ee_n^*\left(x\right)\ee_m^*\left(y\right)\delta_{j,\left(n-1\right)m\left(y\right)+m}.
\end{align*}
Thus, 
\begin{align*}
\supp\left(x\otimes y\right)=&\left\lbrace \left(n-1\right)m\left(y\right)+m: n\in \supp\left(x\right), m\in \supp\left(y\right)\right\rbrace. 
\end{align*}
From this and definition one easily obtains that 
\begin{itemize}
\item $m\left(x\otimes y\right)=m\left(x\right)m\left(y\right)$. 
\item $\left\vert \supp\left(x\otimes y\right)\right\vert= \left\vert \supp\left(x\right)\right\vert \left\vert \supp\left(y\right)\right\vert$.
\item If $x_1, x_2, y_1, y_2\in \spn\left(\EB\right)\setminus \left\lbrace 0\right\rbrace$ and 
\begin{align*}
&m\left(y_1\right)=m\left(y_2\right),
&&\supp\left(x_1\right)\subset\supp\left(x_2\right),&&& \supp\left(y_1\right)\subset\supp\left(y_2\right), 
\end{align*}
then $\supp\left(x_1\otimes y_1\right)\subset \supp\left(x_2\otimes y_2\right)$. 
\item If $m\left(x\right)=m\left(y\right)$ and $\supp\left(y\right)\subset \supp\left(x\right)$, then for each $j\in \NN$, $m\left(x^j\right)=m\left(y^j\right)=\left(m\left(x\right)\right)^j$, and $\supp\left(y^j\right)\subset \supp\left(x^j\right)$. 
\end{itemize}
\end{remark}
\begin{remark}\rm \label{remarksigns2}If $x, y\in \spn\left(\EB\right)$ and there are $\alpha, \beta\in \SSS_{\FF}^{\NN}$ such that $z=\alpha \cdot x$ and $w=\beta \cdot y$, then there is $\gamma\in \SSS_{\FF}^{\NN}$ such that 
\begin{align*}
z\otimes w=&\gamma\cdot \left(x\otimes y\right). 
\end{align*}
Indeed, if $x\not=0\not=y$, one can for example take 
\begin{align*}
\gamma_j:=&
\begin{cases}
\alpha_{n+1}\beta_m & \text{ if } 0\le n< m\left(x\right), 1\le m\le m\left(y\right) \text{ and } j=m\left(y\right)n+m;\\
1 & \text{in any other case}. 
\end{cases}
\end{align*}
\end{remark}

\section{Some basic properties}\label{sectionbasic}

In this section, we state some elementary properties of minimal systems satisfying \eqref{p**} - or weaker, related properties -, and prove some auxiliary results.

\begin{remark}\rm \label{remarklambdal}If $\EB$ satisfies \eqref{p*l-}, then iterating we obtain 
\begin{align*}
&\frac{\left\Vert x\right\Vert^{2^k}}{K^{2^k-1}}\le \left\Vert x^{2^k}\right\Vert &&\forall x\in \spn\left(\EB\right), \forall k\in \NN_0. 
\end{align*}
In particular, this implies that $\lambda\left(1\right)\le K$ and that if there is $n\in \NN$ such that $\lambda\left(n\right)>K$, then $\lambda$ is unbounded. 
\end{remark}

\begin{remark}\rm \label{remarklambdar}If $\EB$ satisfies \eqref{p*-}, then iterating we obtain 
\begin{align*}
\left\Vert x^{2^k}\right\Vert\le K^{2^k-1} \left\Vert x\right\Vert^{2^k}   &&\forall x\in \spn\left(\EB\right), \forall k\in \NN_0. 
\end{align*}
Iterating, and taking into account that $\ee_1^*$ is bounded,  this implies that if $\left\vert \ee_1^*\left(x\right)\right\vert\ge 1$, then $\left\Vert x \right\Vert\ge \frac{1}{K}$. In particular, 
\begin{align*}
&\lambda\ge \frac{1}{K}&&\text{ and } &&&\left\Vert \ee_1^*\right\Vert\le K. 
\end{align*}
\end{remark}

\begin{remark}\label{remarkmax}Let $K_1,K_2>0$, let $\left\Vert\cdot \right\Vert_{\triangleleft}$ and $\left\Vert\cdot \right\Vert_{\circ}$ be a norm and a seminorm on $\mathtt{c}_{00}$ respectively, let $\XX$ be the completion of $\mathtt{c}_{00}$ with the norm 
\begin{align*}
\left\Vert x\right\Vert :=&\max\left\lbrace \left\Vert x \right\Vert_{\triangleleft}, \left\Vert x \right\Vert_{\circ}\right\rbrace &&\forall x\in \mathtt{c}_{00},
\end{align*}
and let $\EB$ be the canonical vector system of $\XX$.\\ Suppose that $\EB$ satisfies one of the following conditions with respect to $\left\Vert\cdot \right\Vert_{\triangleleft}$ and with respect to $\left\Vert\cdot \right\Vert_{\circ}$: \eqref{p*}, \eqref{p*-+}, \eqref{p*-}, \eqref{p*l-+}, \eqref{p*l-}, and for $K=K_1$ and $K=K_2$ respectively.  \\
Then $\EB$ satisfies the same condition with respect to $\left\Vert \cdot\right\Vert$, for $K=\max\left\lbrace K_1,K_2\right\rbrace$. \\
On the other hand, if $\EB$ satisfies \eqref{p**} with respect to both $\left\Vert\cdot \right\Vert_{\triangleleft}$ and $\left\Vert\cdot \right\Vert_{\circ}$ and for some $K=K_0$, then $\EB$ does \emph{not} satisfy \eqref{p*l} with respect to $\left\Vert \cdot\right\Vert$ for any $K>0$, unless $\left\Vert \cdot\right\Vert_{\circ}\lesssim \left\Vert \cdot\right\Vert_{\triangleleft}$ (or $\left\Vert \cdot\right\Vert_{\circ}\gtrsim \left\Vert \cdot\right\Vert_{\triangleleft}$, if $\left\Vert \cdot\right\Vert_{\circ}$ is also a norm). Indeed, if $\EB$ does satisfy \eqref{p*l}, then for all $x, y\in \spn\left(\EB\right)$, 
\begin{align*}
\left\Vert x\right\Vert \left\Vert y\right\Vert\lesssim \max\left\lbrace \left\Vert x\right\Vert_{\triangleleft} \left\Vert y\right\Vert_{\triangleleft}, \left\Vert x\right\Vert_{\circ} \left\Vert y\right\Vert_{\circ}\right\rbrace, 
\end{align*}
whence the result follows. 
\end{remark}

\begin{lemma}\label{lemmabounded}Let $\EB$ be a minimal system for $\XX$. If $\EB$ satisfies \eqref{p*}, then for each $n\in \NN_{\ge 3}$, 
\begin{align*}
\left\Vert \ee_n\right\Vert\le& K\left\Vert \ee_2\right\Vert  \left\Vert \ee_1\right\Vert. 
\end{align*}
Hence, if $\EB^*$ is bounded, then both $\EB$ and $\EB^*$ are seminormalized. 
 \end{lemma}
 \begin{proof}
 Fix $\epsilon>0$, and let $x_{\epsilon}=\ee_1+\epsilon\ee_{n-1}$. Since
\begin{align*}
\left\Vert \ee_n+ \epsilon \ee_{2n-2}\right\Vert=&\left\Vert \ee_2\otimes x_{\epsilon}\right\Vert\le K \left\Vert \ee_2\right\Vert  \left\Vert x_{\epsilon}\right\Vert, 
\end{align*}
the result is obtained letting $\epsilon\rightarrow 0$. 
\end{proof}

\begin{lemma}\label{lemmaboundedbelow}Let $\EB$ be a minimal system for $\XX$. If {$\EB$} satisfies \eqref{p*l}, then for each $n\in \NN_{\ge 3}$,
\begin{align*}
\left\Vert \ee_n\right\Vert\ge& K^{-1}\left\Vert \ee_2\right\Vert  \left\Vert \ee_1\right\Vert. 
\end{align*}
 \end{lemma}
 \begin{proof}
This is proven in the same manner as Lemma~\ref{lemmabounded}.  
 \end{proof}

 \begin{corollary}\label{corollaryseminormalized}Let $\EB$ be a minimal system for $\XX$. If $\EB$ satisfies \eqref{p**}, it is seminormalized. 
 \end{corollary}

\begin{lemma}\label{lemman1/p}Let $\EB$ be a minimal system for $\XX$. If $\EB$ satisfies \eqref{p**}, there is $a\in \left[0, 1\right]$ such that $\lambda\left(n\right)\approx n^{a}$.
\end{lemma}
\begin{proof}
By the proof of \cite[Proposition 2.a.9]{LT1977} (cf. \cite[Proposition 3.2]{CDKM2023}), there are $C>0$ and $a\in \RR$ such that
\begin{align*}
&\frac{n^{a}}{C}\le \lambda\left(n\right)\le C n^a &&\forall n\in \NN. 
\end{align*}
As $\ee_1^*\left(\Ind_{I_n}\right)=1$ for each $n\in \NN$, it follows that $a\ge 0$, whereas Corollary~\ref{corollaryseminormalized} implies that $a\le 1$. 
\end{proof}

\section{Some auxiliary results on translations and invariance under spreading}\label{sectiontranslations}
In this section, we prove some technical results that will be useful in the sequel. {Recall the notation and the terminology introduced in
Subsection~\ref{subsectiontranslations}.}
\begin{remark}\rm \label{remarktranslateseminorm}If a minimal system $\EB$ is $K$-translation invariant, the operator $T_k:\spn\left(\EB\right)\rightarrow \spn\left(\EB_{>k}\right)$ is an isomorphism with $$\max\left\lbrace \left\Vert T_k\right\Vert, \left\Vert T_k^{-1}\right\Vert\right\rbrace\le K.$$ In particular, $\EB$ is seminormalized. 
\end{remark}

\begin{remark}\rm\label{remarkdualsystem} If a minimal system $\EB$ is $K$-right-translation invariant, then $\EB^*$ is $K$-left translation invariant. 
\end{remark}

\begin{remark}\rm\label{remarkdualbasistranslate} If a Schauder basis $\EB$ is $K$-translation invariant, then $\EB^*$ is $K$-left translation invariant and $KK_{b,c}$-right-translation invariant. 
\end{remark}

\begin{example}\label{examplenotbounded1}Let $\XX$ be the completion of $\mathtt{c}_{00}$ with the norm
\begin{align*}
\left\Vert \left(a_n\right)_{n\in\NN}\right\Vert:=&\sup_{n\in\NN}\left\vert a_{n-1}-a_n\right\vert,
\end{align*}
where $a_0:=0$, and let $\EB$ be the canonical vector system of $\XX$. Then $\EB$ is a normalized, $1$-translation invariant minimal system. On the other hand, $\EB^*$ is not bounded, so in particular it is not right-translation invariant.
\end{example}

\begin{lemma}\label{lemmatranslates1} Let $\EB$ be a translation invariant minimal system for $\left(\XX, \|\cdot\|\right)$, let $\alpha:=\limsup_{k\to \infty}\left\Vert \ee_k\right\Vert$, and let $K$ be a translation invariance constant for $\EB$. There is a renorming $\left(\XX, \|\cdot\|_t\right)$ for which $\EB$ is normalized and $1$-translation invariant, with equivalence constant depending only on $K$ and $\alpha$. 
\begin{enumerate}[\rm (i)]
\item \label{dualbounded} If $\EB^*$ is bounded, we can choose the renorming so that $\EB^*$ is normalized, and the equivalence constant depends only on $K$, $\alpha$ and the supremum of the norms of the $\ee_n^*$'s. 
\item \label{schauderstranslate} If $\EB$ is a Schauder basis, we can choose the renorming so that $\EB$ is also bimonotone, and the equivalence constant depends only on $K$, $\alpha$ and $K_b$. 
\item \label{uncondtranslates} If $\EB$ is unconditional, we can choose the renorming so that $\EB$ is $1$-unconditional, and the equivalence constant depends only on $K$, $\alpha$ and $K_u$.
\end{enumerate}
\end{lemma}
\begin{proof}
Define on $\spn\left(\EB\right)$ the norm 
\begin{align}
\left\Vert x\right\Vert_{t}:=&\frac{\limsup_{k\to \infty}\left\Vert T_k\left(x\right)\right\Vert}{\alpha}. \label{translationinvariant}
\end{align}
Since $\EB$ is seminormalized and $K$-translation invariant, $\left\Vert  \cdot \right\Vert_{t}$ is equivalent to  $\left\Vert  \cdot \right\Vert$, so we can extend the former to an equivalent norm on $\XX$, and it is clear that $\EB$ is normalized and $1$-translation invariant with respect to $\left\Vert  \cdot \right\Vert_{t}$.   \\
\ref{dualbounded} First we apply the above procedure, and then we define
\begin{align*}
&\left\Vert x\right\Vert_{\triangleleft}=\max\left\lbrace \left\Vert x\right\Vert_{{\ell_{\infty}}}, \left\Vert x\right\Vert_t\right\rbrace&&\forall x\in \XX. 
\end{align*}
Since $\EB^*$ and $\EB$ are both seminormalized, $\left\Vert x\right\Vert_{\triangleleft}$ is an equivalent norm, and $\EB$ and $\EB^*$ are both normalized with respect to $\left\Vert \cdot \right\Vert_{\triangleleft}$. As $\EB$ is $1$-translation invariant with respect to both $ \left\Vert x\right\Vert_{{\ell_{\infty}}}$ and $ \left\Vert x\right\Vert_{{t}}$, the same holds for $ \left\Vert x\right\Vert_{\triangleleft}$. \\
\ref{schauderstranslate} First we pass to the renorming given by \eqref{translationinvariant} so that $\EB$ is normalized and $1$-translation invariant, and then we use a standard renorming with respect to which $\EB$ is bimonotone. Note that the latter renorming preserves normalization and $1$-translation invariant, so we are done. \\
\ref{uncondtranslates} This is proved as \ref{schauderstranslate}, but with the standard renorming that makes $\EB$ a $1$-unconditional basis. 
 \end{proof}

The next lemma gives us a very weak form of bidemocracy for intervals. Under further conditions, it will allow us to replace bidemocracy by generally weaker hypotheses. 

\begin{lemma}\label{lemmabidemlite}Let $\EB$ be a $K$-translation invariant Schauder basis for $\XX$. The following hold: 
\begin{enumerate}[\rm (i)]
\item \label{projection-light} 
For every $x\in \spn\left(\EB\right)$, every $m\ge m\left(x\right)$ and every $l\in \NN_0$, 
\begin{align*}
\left\vert \sum_{n=1}^{m}\ee_n^*\left(x\right)\right\vert\lVert \Ind_{I_{l+1,l+m}}\rVert \le& C_0 m\left\Vert x\right\Vert, 
\end{align*}
where $C_0=2K\max\left\lbrace K_b,K_{b,c}K\right\rbrace$.
\item \label{bidem-light}
For each $m\in \NN$ and each $l\in\NN_0$, 
\begin{align*}
\lVert \Ind_{I_{l+1, l+m}}\rVert \lVert \Ind_{I_{l+1, l+m}}^*\rVert \le C_1m, 
\end{align*}
where $C_1=C_0KK_bK_{b,c}$. 
\end{enumerate}
\end{lemma}
\begin{proof}
\ref{projection-light} As
\begin{align*}
&\left( \sum_{n=1}^{m}\ee_n^*\left(x\right)\right)  \Ind_{I_{l+1, l+m}}= \sum_{k=0}^{m-1}\left(S_{l+m}\left(T_{k+l}\left(x\right)\right)+T_{m}^{-1}\left(S_{l+m,c}\left(T_{k+l}\left(x\right)\right) \right)\right),
\end{align*}
the triangle inequality gives the desired estimate. \\
\ref{bidem-light}
Given $\epsilon>0$, choose $x\in \SSS_{\XX}$ so that 
\begin{align*}
\left\Vert \Ind_{I_{l+1, l+m}}^*\right\Vert\le&\epsilon +\left\vert \Ind_{I_{l+1, l+m}}^*\left(x\right)\right\vert.
\end{align*}
By \ref{projection-light}, 
\begin{align*}
\left\vert \Ind_{I_{l+1, l+m}}^*\left(x\right)\right\vert=&\left\vert \sum_{n=l+1}^{l+m}\ee_n^*\left(x\right)\right\vert=\left\vert \sum_{n=1}^{m}\ee_n^*\left(T_{l}^{-1}\left(S_{l+m}\left(S_{l,c}\left( x\right)\right)\right)\right)\right\vert\\
\le& \frac{C_0m}{\left\Vert \Ind_{I_{l+1, l+m}}\right\Vert}KK_bK_{b,c},
\end{align*}
and we complete the proof by letting $\epsilon\rightarrow 0$. 
\end{proof}
\begin{remark}\rm \label{remarksimplerbidemlight} Note that using Lemma~\ref{lemmatranslates1} and considering an equivalent norm, one can obtain $C_0=C_1=2$ in Lemma~\ref{lemmabidemlite}. 
\end{remark}

We also want to replace democracy with a generally weaker hypothesis. The next lemma takes care of that. 

\begin{lemma}\label{lemmademoc0}Let $\EB$ be a minimal system for $\XX$ that is right-translation invariant. The following are equivalent: 
\begin{enumerate}[\rm i)]
\item \label{dem1}$\EB$ is democratic. 
\item \label{cons1}$\EB$ is conservative.
\item \label{pardem1}$\EB$ is partially democratic. 
\end{enumerate}
If $\EB$ is also a translation invariant Schauder basis, the following are equivalent: 
\begin{enumerate}[\rm a)]
\item \label{dem2}$\EB^*$ is democratic. 
\item \label{cons2}$\EB^*$ is conservative.
\item \label{pardem2}$\EB^*$ is partially democratic. 
\end{enumerate}
\end{lemma}
\begin{proof}
Suppose $\EB$ is partially democratic with constant $C$ and is right-translation invariant with constant $K$, and choose $A, B\in \NN^{<\infty}$ with $\left\vert A\right\vert\le \left\vert B\right\vert$. By partial democracy, there is $A\subset D\in \NN^{<\infty}$ such that, for each $E\subset \NN\setminus D$ with $\left\vert A\right\vert\le \left\vert E\right\vert<\infty$, $\left\Vert \Ind_A\right\Vert \le C\left\Vert \Ind_E\right\Vert$. Pick $m>\max\left(D\right)$. We have
\begin{align*}
\left\Vert \Ind_A\right\Vert \le C\left\Vert T_{m}\left(\Ind_B\right)\right\Vert\le CK \left\Vert \Ind_B\right\Vert. 
\end{align*}
If $\EB$ is a Schauder basis, the result follows by the above argument and Remark~\ref{remarkdualbasistranslate}. 
\end{proof}

\begin{corollary}\label{corollarybidemtranslate} Let $\EB$ be a Schauder basis for $\XX$. If $\EB$ is translation invariant and both $\EB$ and $\EB^*$ are partially democratic, then $\EB$ is bidemocratic. 
\end{corollary}
\begin{proof}
This follows at once from Lemmas~\ref{lemmabidemlite} and~\ref{lemmademoc0}. 
\end{proof}
Next, we obtain some results on the relations between translation invariance and invariance under spreading (see Definition~\ref{weaklysubsym}).

\begin{lemma}\label{lemmatranslationtosubsym}Let $\EB$ be a minimal system for $\XX$. The following hold: 
\begin{enumerate}[\rm (i)]
\item \label{uniformti}If every subsequence of $\EB$ is translation invariant, there is $K>0$ such that every subsequence of $\EB$ is $K$-translation invariant. 
\item \label{tinvweaksub} If  every subsequence of $\EB$ is $K$-translation invariant, then every two subsequences of $\EB$ are $K^2$-equivalent. Hence, in particular, $\EB$ is $K^2$-invariant under spreading. 
\item \label{weaksub1} If every two subsequences of $\EB$ are equivalent, there is $C>0$ such that every two subsequences of $\EB$ are $C$-equivalent. In particular, $\EB$ is  $C$-invariant under spreading. 
\item \label{weaksubtic}  If $\EB$ is $C$-invariant under spreading, every subsequence of $\EB$ is $C^2$-translation invariant, and every two subsequences of $\EB$ are $C^2$-equivalent. 
\end{enumerate}
\end{lemma}
\begin{proof}
\ref{uniformti} 
First, by Lemma~\ref{lemmatranslates1}, we may assume that $\EB$ is normalized and $1$-translation invariant.   \\

Second, note that it is enough to prove that there are $K>0$ and $n\in \NN$
	such that every subsequence of $\EB_{>n}$ is $K$-translation invariant.
	Indeed, suppose that this is the case, and let
	$\BB=\left(\ee_{g\left(i\right)}\right)_{i\in \NN}$ be an arbitrary subsequence
	of $\EB$, with $g:\NN\rightarrow \NN$ strictly increasing. Set
	\begin{align*}
		&\BB^{\left(n\right)}:=\left(T_n\left(\ee_{g\left(i\right)}\right)\right)_{i\in \NN}=\left(\ee_{g\left(i\right)+n}\right)_{i\in \NN},
		&&\ZB:=\overline{\spn\left(\BB\right)},
		&&&\ZB_n:=\overline{\spn\left(\BB^{\left(n\right)}\right)}.
	\end{align*}
	Since $\left(g\left(i\right)+n\right)_{i\in \NN}$ is strictly increasing and
	$g\left(i\right)+n>n$ for every $i\in \NN$, $\BB^{\left(n\right)}$ is a
	subsequence of $\EB_{>n}$, and is thus $K$-translation invariant by hypothesis.
	On the other hand, given
	$x=\sum_{i\in \NN}a_i\ee_{g\left(i\right)}\in \spn\left(\BB\right)$ and
	$k\in \NN_0$, we have
	\begin{align*}
		T_k\left[\BB^{\left(n\right)},\ZB_n\right]\left(T_n\left(x\right)\right)
		=&\sum_{i\in \NN}a_i\ee_{g\left(i+k\right)+n}
		=T_n\left(T_k\left[\BB,\ZB\right]\left(x\right)\right),
	\end{align*}
	that is, $T_n$ intertwines the translation operators of $\BB$ with those of
	$\BB^{\left(n\right)}$. As $\EB$ is $1$-translation invariant, $T_n$ is an
	isometry on $\spn\left(\EB\right)$, so the above gives
	\begin{align*}
		&\left\Vert x\right\Vert=\left\Vert T_n\left(x\right)\right\Vert
		&&\text{and}
		&&&\left\Vert T_k\left[\BB,\ZB\right]\left(x\right)\right\Vert
		=\left\Vert T_k\left[\BB^{\left(n\right)},\ZB_n\right]\left(T_n\left(x\right)\right)\right\Vert.
	\end{align*}
	Combining the two displays with the $K$-translation invariance of
	$\BB^{\left(n\right)}$ we obtain
	\begin{align*}
		&\frac{1}{K}\left\Vert x\right\Vert\le \left\Vert T_k\left[\BB,\ZB\right]\left(x\right)\right\Vert\le K\left\Vert x\right\Vert
		&&\forall x\in \spn\left(\BB\right),\ \forall k\in \NN_0,
	\end{align*}
	so $\BB$ is $K$-translation invariant, as claimed. 
	
	So, suppose for a
	contradiction that this does not hold, and pick a subsequence of $\EB$ that
	fails to be $1$-translation invariant, say $\EB_1:=\left(\ee_{g_1\left(n\right)}\right)_{n\in \NN}$ with $g_1: \NN\rightarrow \NN$ strictly increasing. Let $\XX_1$ be the closed linear span of $\EB_1$. By our choice, there are 
$y_1\in  \spn\left(\EB_1\right)\cap \SSS_{\XX}$ and $k_1\in \NN$ such that 
\begin{align*}
\max\left\lbrace \frac{1}{\left\Vert T_{k_1}\left[\EB_1,\XX_1\right]\left(y_1\right)\right\Vert}, \left\Vert T_{k_1}\left[\EB_1,\XX_1\right]\left(y_1\right)\right\Vert \right\rbrace>1. 
\end{align*}
Now choose $m_1>m\left[\EB,\XX\right]\left(T_{k_1}\left[\EB_1,\XX_1\right]\left(y_1\right)\right)+1$ (note that the translation index $k_1$ is relative to the subsequence $\EB_1$, but $m\left[\EB, \XX\right]$ is the maximum relative to the original system $\EB$).\\
By hypothesis, there is a 
subsequence of $\EB_{>m_1}$ that fails to be $2$-translation invariant, say $\EB_2:=\left(\ee_{g_2\left(n\right)}\right)_{n\in \NN}$. Using the same notation as before, we get that there is 
$y_2\in  \spn\left(\EB_2\right)\cap \SSS_{\XX}$ and $k_2\in \NN$ such that 
\begin{align*}
\max\left\lbrace \frac{1}{\left\Vert T_{k_2}\left[\EB_2,\XX_2\right]\left(y_2\right)\right\Vert}, \left\Vert T_{k_2}\left[\EB_2,\XX_2\right]\left(y_2\right)\right\Vert \right\rbrace>2. 
\end{align*}
Inductively, our construction yields for each $j$ a subsequence $\left(\EB_j\right)_{j\in \NN}$,  a vector $y_j\in \SSS_{\XX}\cap \XX_j$ and a positive integer $k_j$ such that 
\begin{align*}
\max\left\lbrace \frac{1}{\left\Vert T_{k_j}\left[\EB_j,\XX_j\right]\left(y_j\right)\right\Vert}, \left\Vert T_{k_j}\left[\EB_j,\XX_j\right]\left(y_j\right)\right\Vert \right\rbrace>j. 
\end{align*}
Moreover, for each $j\in \NN$, 
\begin{align*}
m\left[\EB, \XX\right]\left(T_{k_j}\left[\EB_j,\XX_j\right]\left(y_j\right)\right)<g_{j+1}\left(1\right)\le l\left[\EB, \XX\right]\left(y_{j+1}\right).
\end{align*}
Let $\YB=\left(\yy_{i}\right)_{i\in I}$ be the subsequence of $\EB$ constructed by taking 
\begin{align*}
I:=&\bigcup_{j\in \NN}\left\lbrace g_j\left(n\right): 1\le n\le m \left[\EB_j,\XX_j\right]\left(T_{k_j}\left(\EB_j,\XX_j\right]\left(y_j\right)\right)\right\rbrace
\end{align*}
in its natural ordering, and let $\YY:=\overline{\spn\left(\YB\right)}$.  Note that for each $j\in \NN$, 
\begin{align*}
T_{k_j}\left[\YB, \YY\right]\left(y_j\right)=&T_{k_j}\left[\EB_j, \XX_j\right]\left(y_j\right). 
\end{align*}
As $\YB$ is translation invariant, this is a contradiction, and \ref{uniformti} is proved.  \\
\ref{tinvweaksub} Fix $m\in \NN$, let $\left(k_j\right)_{1\le j\le m}\subset \NN$ and $\left(l_j\right)_{1\le j\le m}\subset\NN$ be strictly increasing finite sequences, choose $m_1>\max\left\lbrace k_m,l_m\right\rbrace$, and set $l_j:=k_j:=m_1+j$ for all $m+1\le j\le 2m$. Given $\left(a_j\right)_{1\le j\le m}\subset \FF$, by the $K$-translation invariance of every subsequence of $\EB$ we have
\begin{align*}
\left\Vert \sum_{j=1}^{m}a_j \ee_{k_j}\right\Vert\le K \left\Vert \sum_{j=1}^{m}a_j \ee_{k_{j+m}}\right\Vert= K \left\Vert \sum_{j=1}^{m}a_j \ee_{l_{j+m}}\right\Vert\le K^2\left\Vert \sum_{j=1}^{m}a_j \ee_{l_j}\right\Vert.
\end{align*}
\ref{weaksub1} This is proved by a simplified version of the argument of \ref{uniformti} (alternatively, this can be derived from \ref{uniformti} and~\ref{tinvweaksub}).  \\
\ref{weaksubtic} Clear. 
\end{proof}

\section{From minimal systems to Schauder bases}\label{sectionschauder}

In this section, we find sufficient conditions under which minimal systems satisfying \eqref{p**} are Schauder bases. We begin by showing that \eqref{p**} entails invariance under translations. 

\begin{lemma}\label{lemmatranslates} Let $\EB$ be a minimal system for $\XX$. If $\EB$  satisfies  \eqref{p**}, then $\EB$ is translation invariant, and $\EB^*$ is left-translation invariant. More precisely,  for each $x\in \spn\left(\EB\right)$ and each $k, m\in \NN_0$, 
\begin{align}
\left\Vert T_k\left[\EB,\XX\right]\left(x\right)\right\Vert \le K^2 \left\Vert T_m\left[\EB,\XX\right]\left(x\right)\right\Vert,  \label{translates2} 
\end{align}
whereas for each $x^*\in \spn\left(\EB^*\right)$ and each $k\le m\in \NN_0$,
\begin{align}
\left\Vert T_k\left[\EB^*,\XX^*\right]\left(x^*\right)\right\Vert \le K^2 \left\Vert T_m\left[\EB^*,\XX^*\right]\left(x^*\right)\right\Vert.  \label{translates2.5}
\end{align}
Moreover,$K^{-1}\le \left\Vert \ee_n\right\Vert\le K$ and 
\begin{align}
&\frac{1}{K^3}\le \left\Vert \ee_n\right\Vert\le K^3 &&\forall n\in \NN_{\ge 2}.  \label{seminormalized2}
\end{align}
\end{lemma}
\begin{proof}
 Pick $k,m\in \NN_0$ and $0\not=x\in \spn\left(\EB\right)$, let $j_0:=1+m\left(x\right)$, choose $\epsilon>0$, and let
$y_{\epsilon}:=x+\epsilon \ee_{j_0+2m}$ and $z_{\epsilon}:=x+\epsilon \ee_{j_0+2k}$. By Remark~\ref{remarktranslate}, 
\begin{align*}
&\ee_2\otimes T_k\left(y_{\epsilon}\right)=T_{j_0+2m+k}\left(T_k\left(y_{\epsilon}\right)\right)=T_{j_0+2m+2k}\left(y_{\epsilon}\right);\\
&\ee_2\otimes T_m\left(z_{\epsilon}\right)=T_{j_0+2k+m}\left(T_m\left(z_{\epsilon}\right)\right)=T_{j_0+2m+2k}\left(z_{\epsilon}\right). 
\end{align*}
Hence, 
\begin{align*}
\left\Vert T_k\left(y_{\epsilon}\right)\right\Vert\le& \frac{K}{\left\Vert \ee_2\right\Vert}\left\Vert \ee_2\otimes T_k\left(y_{\epsilon}\right)\right\Vert= \frac{K}{\left\Vert \ee_2\right\Vert}\left\Vert T_{j_0+2m+2k}\left(y_{\epsilon}\right)\right\Vert\\
=& \frac{K}{\left\Vert \ee_2\right\Vert}\frac{\left\Vert T_{j_0+2m+2k}\left(y_{\epsilon}\right)\right\Vert}{\left\Vert T_{j_0+2m+2k}\left(z_{\epsilon}\right)\right\Vert}\left\Vert \ee_2\otimes T_m\left(z_{\epsilon}\right)\right\Vert \\
\le& K^2\frac{\left\Vert T_{j_0+2m+2k}\left(y_{\epsilon}\right)\right\Vert}{\left\Vert T_{j_0+2m+2k}\left(z_{\epsilon}\right)\right\Vert}\left\Vert T_m\left(z_{\epsilon}\right)\right\Vert.
\end{align*}
Letting $\epsilon\rightarrow 0$ we obtain \eqref{translates2}, whereas \eqref{translates2.5} follows from this and Remark~\ref{remarkdualsystem}.\\
Finally, the upper and lower bounds for the norm of $\ee_1$ follow from Remarks~\ref{remarklambdal} and~\ref{remarklambdar}, whereas \eqref{seminormalized2} follows from said bounds and \eqref{translates2}. 
\end{proof}
\begin{corollary}\label{corollaryp**1tinv}Let $\EB$ be a minimal system for $\XX$ satisfying \eqref{p**}. There is an equivalent norm on $\XX$ for which $\EB$ is normalized and $1$-translation invariant, and the equivalence constant depends only on $K$. 
\end{corollary}
\begin{proof}
This follows from Lemmas~\ref{lemmatranslates1} and~\ref{lemmatranslates}. \end{proof}

We continue with a relation that will be useful  in the sequel. As it follows easily from definitions, the proof is left to the reader. 
\begin{lemma}\label{lemmamovingright}Let $\EB$ be a minimal system for $\XX$, and $x,y, z\in \spn\left(\EB\right)\setminus \left\lbrace 0\right\rbrace$. If $m\left(x\right)=m\left(y\right)=:m$, then 
\begin{align*}
z\otimes \left(x+T_m\left(y\right)\right)=&\ee_1^*\left(z\right) x+T_{m}\left( \left(T_{1}^{-1}\left(S_{1,c}\left(z\right)\right) \otimes \left(y+T_m\left(x\right)\right)\right)\right)\\
&+T_m\left(\sum_{j=1}^{m\left(z\right)}\left(\ee_j^*\left(z\right)-\ee_{j+1}^*\left(z\right)\right)T_{2m\left(j-1\right)}\left(y\right)\right)\\
=&\ee_1^*\left(z\right) x+  \left(z - T_{1}^{-1}\left(S_{1,c}\left(z\right)\right)\right)  \otimes T_m\left(y\right)\\
&+T_{m}\left( \left(T_{1}^{-1}\left(S_{1,c}\left(z\right)\right) \otimes \left(y+T_m\left(x\right)\right)\right)\right). 
\end{align*}
\end{lemma}
Now we introduce an auxiliary norm. 

\begin{definition}\label{definitiondifferences}Let $\EB$ be a minimal system for $\XX$. We define $\cD:\spn\left(\EB\right)\rightarrow \spn\left(\EB\right)$  by 
\begin{align*}
\cD\left(x\right):=&x-T_1^{-1}\left(S_{1,c}\left(x\right)\right)=\sum_{n=1}^{m\left(x\right)}\left(\ee_n^*\left(x\right)-\ee_{n+1}^*\left(x\right)\right)\ee_n,
\end{align*}
and the norm $N_{\cD}=N_{\cD}\left[\EB, \XX, \left\Vert \cdot\right\Vert\right]$ by 
\begin{align*}
N_{\cD}\left(x\right):=&\left\vert \ee_1^*\left(x\right)\right\vert+ \left\Vert \cD\left(x\right) \right\Vert. 
\end{align*}
\end{definition}

\begin{remark}\label{remarkNDsmaller}
If a minimal system is translation invariant, $N_{\cD}$ is dominated by the norm of the space. More precisely, if $\EB$ is $K$-translation  invariant, the triangle inequality gives  
\begin{align}
N_{\cD}\left(x\right)\le& \left(\left\Vert \ee_1^*\right\Vert + 1+K +K\left\Vert \ee_1^*\right\Vert \left\Vert \ee_1\right\Vert \right)\left\Vert x\right\Vert \nonumber
\end{align}
for all $x\in \spn\left(\EB\right)$. 
\end{remark}
\begin{remark}\label{remarklambdabounded}Note that for every $n\in \NN$, 
\begin{align*}
N_{\cD}\left(\Ind_{I_{n}}\right)=1+\left\Vert  \ee_n\right\Vert. 
\end{align*}
\end{remark}
 
\begin{lemma}\label{lemmaNDsmaller}
Let $\EB$ be a minimal system for $\XX$. If $\EB$ is translation invariant and $N_{\cD}$ is an equivalent norm, then there is $C\ge 1$ such that
\begin{align*}
&\left\Vert \ee_n^*\right\Vert\ge \frac{n^{k}}{C^k k!}&&\forall n,k\in \NN.
\end{align*}
\end{lemma}
\begin{proof}
By Lemma~\ref{lemmatranslates1}, we may  assume that $\EB$ is normalized and $1$-translation invariant. Choose $C_0\ge 1$ so that 
\begin{align*}
C_0^{-1}N_{\cD}\left(x\right)\le   \left\Vert x\right\Vert\le& C_0 N_{\cD}\left(x\right)
\end{align*}
for all $x\in \spn\left(\EB\right)$. By Remark~\ref{remarklambdabounded}, $\lambda\left(n\right)\le 2C_0$ for each $n\in \NN$. \\
For every $0\not=\alpha \in \NN_{0}^{\NN}$ and every $n\in  \NN$, set 
\begin{align*}
x_{\alpha,n}:=& \sum_{k=1}^{n}\alpha_k\ee_k;\\
y_{\alpha,n}:=&\sum_{k=1}^{n}s_k\left(\alpha\right)\ee_k+\sum_{k=n+1}^{n+ s_n\left(\alpha\right)-1}\left(s_n\left(\alpha\right)-k+n\right)\ee_k. 
\end{align*} 
Given $\alpha, n$ as above, 
\begin{align*}
\left\Vert y_{\alpha,n}\right\Vert\le& C_0 N_{\cD}\left(y_{\alpha,n}\right)=C_0\left\vert \alpha_1\right\vert+C_0\left\Vert -\sum_{k=1}^{n-1}\alpha_{k+1}\ee_k+\sum_{k=n}^{n+ s_n\left(\alpha\right)-1}\ee_k\right\Vert\\
\le&  C_0\left\vert \alpha_1\right\vert+C_0\left\Vert T_1^{-1}S_{1,c}\left(x_{\alpha,n}\right)\right\Vert +C_0 \lambda\left( s_n\left(\alpha\right)\right)\\
\le& 2C_0\left\vert \alpha_1\right \vert +C_0\left\Vert x_{\alpha,n}\right\Vert+2C_0^2 \le  \left(2C_0\left\Vert \ee_1^*\right\Vert+C_0\right)\left\Vert x_{\alpha,n}\right\Vert+2C_0^2.
\end{align*}
Let $C_1:=2\left(2C_0\left\Vert \ee_1^*\right\Vert+C_0+2C_0^2+1\right)$.
The above estimate gives 
\begin{align}
\left\Vert y_{\alpha, n}\right\Vert\le& C_1\max\left\lbrace  \left\Vert x_{\alpha,n}\right\Vert, 1\right\rbrace.  \label{boundforyalpha}
\end{align}
for all $0\not=\alpha\in \NN_{0}^{\NN}$ and all $n\in \NN$.\\
Given $0\not=f_0\in \NN_{0}^{\NN}$ finitely supported, let $n_0:=\max\left(\supp\left(f_0\right)\right)$, and recursively construct a sequence $\left(f_j\right)_{j\in \NN_0}$ of finitely supported non-zero elements in $\NN_{0}^{\NN}$ as follows: If $\left(f_j\right)_{0\le j\le j_0-1}$ have been defined for some $j_0\in \NN$ and $n_j=\max\left(\supp\left(f_{j}\right)\right)$ for every $0\le j\le j_0-1$, define $f_{j_0}\in \NN_{ 0}^{\NN}$ by 
\begin{align*}
f_{j_0}\left(l\right)=&\ee_l^*\left(y_{f_{j_0-1}, n_{j_0-1}}\right) &&\forall l\in \NN,
\end{align*}
and let $n_{j_0}:=\max\left(\supp\left(f_{j_0}\right)\right)$. It follows by construction that for each $j\in \NN$, 
\begin{align*}
x_{f_j,n_j}=&y_{f_{j-1}, n_{j-1}}. 
\end{align*}
From the above and \eqref{boundforyalpha}, inductively we obtain 
\begin{align}
\left\Vert x_{f_j,n_j}\right\Vert\le& C_1^{j}\max\left\lbrace \left\Vert x_{f_0,n_0}\right\Vert, 1\right\rbrace &&\forall j\in \NN_0. \label{boundforfk}
\end{align}
On the other hand, for each $1\le n\le n_0$ and each $j\in \NN$, 
\begin{align}
f_j\left(n\right)=&\sum_{k=1}^{n}f_{j-1}\left(k\right). \label{recursion1}
\end{align}
Now we choose $f_0$: pick $n_0\in \NN_{\ge 2}$, and let 
\begin{align*}
f_0\left(n\right):=&
\begin{cases}
1 & \text{ if } 1\le n\le n_0;\\
0 & \text{ if } n>n_0. 
\end{cases}
\end{align*}
 From   \eqref{recursion1}  and the fact that 
\begin{align*}
&\sum_{k=1}^{n}k^j\ge \frac{n^{j+1}}{j+1}&&\forall n, j\in \NN_0, 
\end{align*}
an inductive argument gives 
\begin{align*}
&f_j\left(n\right)\ge \frac{n^j}{j!}&&\forall 1\le n\le n_0, \forall j\in \NN_0. 
\end{align*}
Hence, using \eqref{boundforfk} and the fact that $\lambda\le 2C_0$ we obtain 
\begin{align*}
&\left\Vert \ee_n^*\right\Vert\ge  \left\vert \ee_n^*\left(\frac{x_{f_j, n_j}}{2C_0C_1^j} \right) \right\vert \ge \frac{n^j}{2C_0C_1^j j!}&&\forall 1\le n\le n_0. 
\end{align*}
As $n_0$ is arbitrary, the proof is complete. 
\end{proof}

Next, we prove the main result of this section.

\begin{lemma}\label{lemmaSchauder3}Let $\EB$ be a minimal system for $\XX$ satisfying \eqref{p**}. The following are equivalent: 
\begin{enumerate}[\rm (i)]
\item \label{bigfunctionals} There is $C>0$ such that 
\begin{align*}
&\left\Vert \ee_n^*\right\Vert\ge \frac{n^k}{C^k k!}  &&\forall n, k\in \NN. 
\end{align*}
\item \label{notsmallfunctionals} $\EB^*$ is unbounded. 
\item \label{notschauder}  $\EB$ is not a Schauder basis. 
\item \label{NDequivalent} $N_{\cD}$ is an equivalent norm. 
\end{enumerate}
If any (and thus all) of the above conditions is satisfied, then $\lambda\le K$ and there is $C>0$ such that
\begin{align}
\varphi\left(n\right)\ge & Cn &&\forall n\in \NN. \label{bigphi}
\end{align}
\end{lemma}
\begin{proof}
By Lemma~\ref{lemmatranslates}, $\EB$ is translation invariant. Hence, by Lemma~\ref{lemmatranslates1}, for the proofs of the equivalences and \eqref{bigphi} we may and will pass to an equivalent norm and assume that $\EB$ is normalized and $1$-translation invariant. While this procedure might change the constant $K$ in \eqref{p**}, once we have established the equivalences we can go back to the original norm and apply \ref{NDequivalent} and Remarks~\ref{remarklambdal} and ~\ref{remarklambdabounded}  to conclude that $\lambda$ is bounded above by the original $K$. \\
\ref{notsmallfunctionals} $\Rightarrow$ \ref{notschauder} Since $\EB$ is normalized, if it were a Schauder basis then $\EB^*$ would be seminormalized. \\
\ref{notschauder} $\Rightarrow$ \ref{NDequivalent}  Suppose, for a contradiction, that $N_{\cD}$ is not an equivalent norm, and set $K_0:=K+1$. By Remark~\ref{remarkNDsmaller}, there is  sequence $\left(z_j\right)_{j\in \NN}\subset \SSS_{\XX}\cap \spn\left(\EB\right)$ such that  
\begin{align*}
&N_{\cD}\left(z_j\right)<\frac{1}{5 jK_0^2}&&\forall j\in \NN. 
\end{align*}
As $\EB$ is not a Schauder basis, for each $j\in \NN$ there exist $x_j, y_j\in \spn\left(\EB\right)$ such that 
\begin{align*}
l\left(y_j\right)>m\left(x_j\right)\ge & 1;\\
\left\Vert x_j\right\Vert=&1; \\
\left\Vert x_j+y_j\right\Vert\le& \frac{1}{4jK_0^2\left(1+\left\Vert z_j\right\Vert_{\ell_1}\right)}.
\end{align*}
For each $j\in \NN$, let $m_j:=m\left(y_j\right)+l\left(y_j\right)-1$, pick
\begin{align*}
0<\epsilon_j<  \frac{1}{4jK_0^2\left(1+\left\Vert z_j\right\Vert_{\ell_1}\right)}\min\left\lbrace 1, \left\vert \ee_{m\left(x_j\right)}^*\left(x_j\right)\right\vert\right\rbrace,
\end{align*}
and set 
\begin{align*}
u_{j}:=&T_{m\left(y_j\right)}\left(x_j\right)+\epsilon_j \ee_{m_j};\\
v_{j}:=&T_{l\left(y_j\right)-1}^{-1}\left(y_j\right)+\epsilon_j \ee_{m_j}. 
\end{align*}
Note that 
\begin{align}
m\left(u_{j}\right)=m\left(v_{j}\right)=m\left(u_{j}+v_{j}\right)=m_j,\label{forproductschauder}
\end{align}
and 
\begin{align*}
u_{j}+T_{m_j}\left(v_{j}\right)=&\epsilon_j\left(\ee_{m_j}+\ee_{2m_j}\right)+T_{m\left(y_j\right)}\left(x_j+y_j\right),
\end{align*}
so 
\begin{align}
\left\Vert u_{j}+T_{m_j}\left(v_{j}\right)\right\Vert\le& 2\epsilon_j +\left\Vert x_j+y_j\right\Vert\le \frac{3}{4jK_0^2\left(1+\left\Vert z_j\right\Vert_{\ell_1}\right)}.\label{firstboundleft}
\end{align}
By Lemma~\ref{lemmamovingright}, the triangle inequality, \eqref{p**}, and construction,  for each $j\in \NN$ we have 
\begin{align*}
\left\Vert v_{j}+T_{m_j}\left(u_{j}\right)\right\Vert\le&  K_0  \left\Vert z_j\otimes\left( v_{j}+T_{m_j}\left(u_{j}\right)\right)\right\Vert \\
\le& K_0 \left\Vert \ee_1^*\left(z_j\right) x_j \right\Vert  +K_0^2\left\Vert z_j-T_1^{-1}\left(S_{1,c}\left(z_j\right)\right)\right\Vert\left\Vert T_{m_j}\left(u_j\right)\right\Vert \\
+&K_0\left\Vert T_{m_j} \left(T_{1}^{-1}\left(S_{1,c}\left(z_j\right)\right)\right) \otimes \left(u_{j}+T_{m_j}\left(v_{j}\right)\right)\right\Vert\\
\le& \left(1+\epsilon_j\right) K_0^2 N_{\cD}\left(z_j\right)+ K_0\left\Vert z_j\right\Vert_{\ell_1} \left\Vert u_{j}+T_{m_j}\left(v_{j}\right)\right\Vert\le\frac{1}{j}. 
\end{align*}
Hence, the triangle inequality gives
\begin{align*}
\left\Vert \ee_2\otimes \left(u_j+v_j\right)\right\Vert=&\left\Vert  u_{j}+T_{m_j}\left(v_{j}\right)+ v_{j}+T_{m_j}\left(u_{j}\right)\right\Vert \le \frac{2}{j}. 
\end{align*}
Thus, 
\begin{align*}
\left\Vert u_j+v_j\right\Vert \le& K_0\left\Vert \ee_2\otimes \left(u_j+v_j\right)\right\Vert \le  \frac{2K_0}{j}. 
\end{align*}
Since 
\begin{align*}
\left\Vert u_{j}\right\Vert\le&  \frac{K_0}{\left\Vert \ee_1-\ee_2\right\Vert}\left\Vert  \left(\ee_1-\ee_2\right)\otimes u_{_j}\right\Vert=  \frac{K_0}{\left\Vert \ee_1-\ee_2\right\Vert}\left\Vert u_{j}-T_{m_j}\left(u_{j}\right)\right\Vert \\
\le& \frac{K_0}{\left\Vert \ee_1-\ee_2\right\Vert}\left\Vert  u_j+T_{m_j}\left(v_{j}\right)\right\Vert + \frac{K_0}{\left\Vert \ee_1-\ee_2\right\Vert}\left\Vert T_{m_j}\left(u_{j}+v_{j}\right)\right\Vert\\
=& \frac{K_0}{\left\Vert \ee_1-\ee_2\right\Vert}\left\Vert  u_j+T_{m_j}\left(v_{j}\right)\right\Vert + \frac{K_0}{\left\Vert \ee_1-\ee_2\right\Vert}\left\Vert u_{j}+v_{j}\right\Vert,
\end{align*}
it follows that 
\begin{align*}
\limsup_{j\to \infty}\left\Vert u_{j}\right\Vert=&0, 
\end{align*}
thus 
\begin{align*}
1=&\limsup_{j\to \infty}\left\Vert x_{j}\right\Vert=\limsup_{j\to \infty}\left\Vert T_{m\left(y_j\right)}\left(x_{j}\right)\right\Vert=\limsup_{j\to\infty}\left\Vert u_j-\epsilon_j\ee_{m_j}\right\Vert=0. 
\end{align*}
\ref{NDequivalent} $\Rightarrow$ \ref{bigfunctionals} This follows by Lemma~\ref{lemmaNDsmaller}. \\
This concludes the proof of the equivalences, and the upper bound for $\lambda$ is proved as stated before. \\ 
It remains to prove \eqref{bigphi}. To this end, we will use the notation and part of the argument of Lemma~\ref{lemmaNDsmaller}: let $\alpha_j:=1$ for all $j\in \NN$ and, for each $n\in \NN_{\ge 2}$, let $x_{\alpha,n}$ and $y_{\alpha, n}$ be as in the proof of Lemma~\ref{lemmaNDsmaller}. Also choose $C_0$ as in that proof. We have 
\begin{align*}
\left\Vert y_{\alpha, n}\right\Vert\le& C_0 N_{\cD}\left(y_{\alpha, n}\right)=C_0\left( 1 +\left\Vert -\Ind_{I_{n-1}}+T_{n-1}\left(\Ind_{I_n}\right)\right\Vert\right)\le 5C_0^2 \;\forall n\in \NN_{\ge 2}. 
\end{align*}
For each $n, m\in \NN_{\ge 2}$ and $1\le k\le m$, let 
\begin{align*}
&A_{m,n,k}:=\left\lbrace mn+ j k : 1-n\le j\le n-1  \right\rbrace,
\end{align*}
and let $B_{m,n,k}:= A_{m,n,k} \setminus\left\lbrace mn\right\rbrace$. Note that 
\begin{align*}
\supp\left(T_{n\left(m-k\right)}\left(y_{\alpha, n}\otimes \ee_k \right) \right)=&A_{m,n,k}. 
\end{align*}
Choose $x_{m,n}^*\in S_{\XX^*}$ so that 
\begin{align*}
x_{m,n}^*\left(\ee_{mn}\right)=& 1. 
\end{align*}
We have 
\begin{align*}
5C_0^2 K\ge& K \left\Vert y_{\alpha, n}\right\Vert\ge   \left\Vert y_{\alpha, n}\otimes \ee_k \right\Vert = \left\Vert T_{n\left(m-k\right)}\left( y_{\alpha, n}\otimes \ee_k \right)\right\Vert\\
\ge& \left\vert x_{m,n}^*\left(T_{n\left(m-k\right)}\left( y_{\alpha, n}\otimes \ee_k \right)\right)\right\vert\\
\ge&  n-\sum_{\substack{1-n\le j\le n-1\\j\not=0}} \left\vert \left(n-\left\vert j\right\vert\right) x_{m,n}^*\left(\ee_{mn+jk} \right) \right\vert \\
\ge& n -n^2\max_{l\in B_{m,n,k}}\left\vert x_{m,n}^*\left(\ee_{l} \right) \right\vert.
\end{align*}
Thus, 
\begin{align}
&\max_{l\in B_{m,n,k}}\left\vert x_{m,n}^*\left(\ee_{l} \right) \right\vert\ge  \frac{1}{2 n }&& \forall n\ge 10C_0^2K, m\ge 2, 1\le k\le m. \label{forcontphinotn1}
\end{align}
Notice that 
\begin{align*}
&B_{m,n,k} \subset \left[ mn+k-kn, mn-k\right] \cup  \left[mn+k, mn+kn-k\right]  \\ &\forall n,m\ge 2,\;\forall  1\le k\le m. 
\end{align*}
Thus, for each $n, i\ge 2$, the family 
\begin{align*}
\left\lbrace B_{n^{i-1} ,n, n^{j-1}} \right\rbrace_{1\le j\le i}
\end{align*}
is pairwise disjoint. \\
Now fix $n_0\ge 10C_0^2K$. For each $i\ge 2$ and each $1\le j\le i$, choose $l_{j,i}\in B_{n_0^{i-1} ,n_0, n_0^{j-1}}$ so that the maximum in \eqref{forcontphinotn1} for $m=n_0^{i-1}$, $n=n_0$ and $k=n_0^{j-1}$ is attained when $l=l_{j,i}$. For each $i\ge 2$, 
\begin{align*}
\left\vert D_i:=\left\lbrace  l_{j,i}\right\rbrace_{1\le j\le i} \right\vert =& i 
\end{align*}
because if $j\not=j'$, then $l_{j,i}$ and $l_{j',i}$ lie in different sets of a pairwise disjoint family. It follows from this that for every $i\ge 2$,
\begin{align*}
4\varphi\left(i\right)\ge& \varphi_u\left(i\right)\ge \max_{\varepsilon\in \FF^{D_i}}\left\Vert \Ind_{\varepsilon, D_i}\right\Vert \ge \max_{\varepsilon\in \FF^{D_i}}\left\vert x^*_{n_0^{i-1}, n_0}\left( \Ind_{\varepsilon, D_i}\right)\right\vert\ge \frac{i}{2 n_0}. 
\end{align*}
\end{proof}

To finish this section, we note that in the proof of Lemma~\ref{lemmaSchauder3} we used both sides of \eqref{p**}. The following example shows that if we relax \eqref{p*l}, the situation is very different, even if we add a translation invariance condition. 

\color{black}

\begin{example}\label{examplenotproduct}Let $\XX$ be the completion of $\mathtt{c}_{00}$ with the norm 
\begin{align*}
\left\Vert \left(a_n\right)_{n\in \NN}\right\Vert=&\max\left\lbrace \left\Vert \left(a_n\right)_{n\in \NN}\right\Vert_{\ell_{\infty}}, \left\vert \sum_{n\in \NN}a_n\right\vert\right\rbrace,
\end{align*}
and let $\EB$ be the canonical vector system of $\XX$. Then $\EB$ is a minimal system with the following properties:
\begin{enumerate}[\rm (1)]
\item \label{(1)} Both $\EB$ and $\EB^*$ are normalized. 
\item \label{(2.5)} $\lambda\left(n\right)=n$ for all $n\in \NN$. 
\item  \label{(2)} $\EB$ is $1$-translation invariant. 
\item  \label{(3)} $\EB$ has properties \eqref{p*} and \eqref{p*l-+} with constant $K=1$.
\item  \label{(4)}$\EB$ is not a Markushevich basis.  
\end{enumerate}
Indeed, \ref{(1)}, \ref{(2.5)} and \ref{(2)} follow at once from definitions; \ref{(3)} follows from Remark~\ref{remarkmax}, whereas \ref{(4)}  follows from the fact that $\left(\frac{\Ind_{I_{n^n+1,n^n+n}}}{n}\right)_{n\in \NN}$ is a normalized Cauchy sequence. 
\end{example}

\color{black}

\section{Some auxiliary results on truncation quasi-greedy minimal systems}\label{sectiontqg}

In this section, we prove some technical results showing that \eqref{p**} (or even less) allow us to go from truncation quasi-greediness to unconditionality.

\begin{lemma}\label{lemmauncondparam2}Let $\EB$ be a minimal system for $\XX$ that satisfies \eqref{p*-} and \eqref{p*l-}. One of the following holds: 
\begin{enumerate}[\rm i)]
\item \label{uncond12}$\EB$ is an unconditional Schauder basis with constant $K_{u}\le K^2$.
\item \label{verycond2} $\EB$ is a conditional minimal system, and there are $0<\delta<1$ and $n\in \NN$ such that
\begin{align*}
&\ku_{n^{2^j}}\ge \left(1+\delta\right)^{2^j}&& \forall j\in \NN. 
\end{align*}
\end{enumerate}
\end{lemma}
\begin{proof}
If \ref{uncond12} does no hold, there is $0<\delta<1$ such that $\EB$ is not unconditional with $K_{u}\le \left(1+\delta\right)K^2$. Thus, there are $x\in \spn\left(\EB\right)\cap \SSS_{\XX}$ and $\alpha\in \SSS_{\FF}^{\NN}$ such that  $\left\Vert y:=\alpha \cdot x\right\Vert > \left(1+\delta\right)K^{2}$. Iterating, we obtain that for every $j\in \NN$, 
\begin{align*}
\left\Vert x^{2^j}\right\Vert\le K^{2^j-1} &&\text{and} &&&\left\Vert y^{2^j}\right\Vert\ge \frac{1}{K^{2^j-1}}\left\Vert y \right\Vert^{2^j}> \frac{\left(1+\delta\right)^{2^j} K^{2^{j+1}}}{K^{2^j-1}},
\end{align*}
so 
\begin{align*}
\left\Vert y ^{2^j}\right\Vert\ge \left(1+\delta\right)^{2^j}K^{2}\left\Vert x^{2^j}\right\Vert.
\end{align*}
Since $y^{2^j}=\alpha^{2^j} \cdot x^{2^j}$, the statement holds for $n=m\left(x\right)$. 
\end{proof}

To apply Lemma~\ref{lemmauncondparam2} to the cases we are interested in, we need the following result from the literature: 
\begin{lemma} \label{lemmalpu1}Let $\EB$ be a minimal system for $\XX$. If $\EB$ is seminormalized and truncation quasi-greedy, then $\kk_m\left[\EB,\XX\right]\lesssim 1+\log m$ and $\EB^*$ is bounded. 
\end{lemma}
\begin{proof}
It follows from \cite[Lemma 2.6]{AAB2023a} that $\EB^*$ is seminormalized, and then the upper bound for $\kk_m$ is given by \cite[Theorem 5.3]{AAW2021}. 
\end{proof}

From the above we obtain the following: 

\begin{corollary}\label{corollaryuncond12}Let $\EB$ be a minimal system for $\XX$ satisfying \eqref{p*-} and \eqref{p*l-}. If $\EB$ is seminormalized and truncation quasi-greedy, then $\EB$ is a $K^2$-unconditional basis. In particular, every quasi-greedy minimal system satisfying \eqref{p**} is $K^2$-unconditional. 
\end{corollary}
\begin{proof}This follows from Lemmas~\ref{lemmauncondparam2} and~\ref{lemmalpu1}, and Corollary~\ref{corollaryseminormalized}. 
\end{proof}

The following boundedness result from the literature will allow us to rule out some of the cases of Lemma~\ref{lemmaSchauder3}.

\begin{lemma}\label{lemmanorming}Let $\EB$ be a seminormalized norming Markushevich basis for $\XX$. If $\EB^*$ is truncation quasi-greedy, it is seminormalized. 
\end{lemma}
\begin{proof}
 By \cite[Lemma 2.6]{AAB2023a},
\begin{align*}
\sup_{n\in \NN}\left\Vert \ee_n^*\right\Vert\left\Vert \ee_n^{**}\big|_{\overline{\spn\left(\EB^*\right)}} \right\Vert<\infty 
\end{align*}
By the norming hypothesis, $\EB$ is equivalent to $\EB^{**}\big|_{\overline{\spn\left(\EB^*\right)}}$, so we are done. 
\end{proof}

\begin{corollary}\label{corollaryuncond2}Let $\EB$ be a seminormalized norming Markushevich basis for $\XX$. If $\EB$ satisfies \eqref{p*-} and \eqref{p*l-} and $\EB^*$ is truncation quasi-greedy, then $\EB$ is $K^2$-unconditional. 
\end{corollary}
\begin{proof}By Lemma~\ref{lemmanorming}, $\EB^*$ is seminormalized, whereas  Lemma~\ref{lemmalpu1} entails $\kk_m\left[\EB^*, \XX^*\right]\lesssim 1+\log m$. Since $\spn\left(\EB^*\right)$ is norming for $\XX$, it follows that $\kk_m\left[\EB, \XX\right]\lesssim 1+\log m$. By Lemma~\ref{lemmauncondparam2}, $\EB$ is $K^2$-unconditional. 
\end{proof}

\section{Main result}\label{mainresults}
In this section, we prove our main result.

\begin{proposition}\label{propositionnewcharac}Let $\EB$ be a minimal system for $\XX$ satisfying \eqref{p**}. Suppose that one of the following conditions hold: 
\begin{enumerate}[\rm (i)]
\item \label{bothpartial}Both $\EB$ and $\EB^*$ are partially democratic. 
\item \label{tqgpartially} $\EB$ is truncation quasi-greedy and partially democratic. 
\item \label{tqgpartially*} $\EB$ is truncation quasi-greedy, and $\EB^*$ is partially democratic. 
\item \label{tqg*partially} $\EB^*$ is truncation quasi-greedy, and $\EB$ is partially democratic. 
\item \label{tqg*dem}$\EB^*$ is truncation quasi-greedy and democratic. 
\item \label{tqg*partially*} $\EB^*$ is truncation quasi-greedy and partially democratic, and one of the following holds: 
\begin{enumerate}[\rm (a)]
\item \label{polynom} For every $c>1$ there exists $n_0$ so that $||\ee_n^*|| \leq c^{n}$ for $n>n_0$.
\item \label{lambdasmall2} $\left\Vert \lambda\right\Vert_{\ell_{\infty}}>K$. 
\item \label{phinotbig} The sequence $\left(\frac{n}{\varphi\left(n\right)}\right)_{n\in \NN}$ is unbounded. 
\item \label{norming2} $\EB$ is a norming Markushevich basis. 
\end{enumerate}
\end{enumerate}
Then $\EB$ is a Schauder basis equivalent to the canonical unit vector basis of $\mathtt{c}_0$ or $\ell_p$ for some $1\le p<\infty$. 
\end{proposition}
\begin{proof}
By Lemma~\ref{lemmatranslates} and Remark~\ref{remarktranslateseminorm}, $\EB$ is seminormalized and translation invariant. \\
Let us consider each of the conditions of the statement: \\
\ref{bothpartial} 
By Lemma~\ref{lemmademoc0}, $\EB$ is democratic. 
Since $\lambda\approx \varphi$, Lemma~\ref{lemmaSchauder3} entails that $\EB$ is a Schauder basis. By Lemma~\ref{lemmatranslates} and Corollary~\ref{corollarybidemtranslate}, $\EB$ is bidemocratic. Now the result follows from \cite[Theorem 3.1]{CDKM2023}. \\
\ref{tqgpartially} By Lemmas~\ref{lemmademoc0} and~\ref{lemmatranslates}, $\EB$ is democratic. By Corollary~\ref{corollaryuncond12}, $\EB$ is unconditional. The result now follows from \cite[Theorem 5.1]{CDKM2023}.\\
\ref{tqgpartially*} By Corollary~\ref{corollaryuncond12}, $\EB$ is unconditional. By Lemmas~\ref{lemmademoc0}, ~\ref{lemmatranslates}, Remark~\ref{remarkdualbasistranslate} and hypothesis, $\EB^*$ is translation invariant and democratic. By Lemma~\ref{lemman1/p}, there is $1\le p\le \infty$ such that $\lambda\left(n\right)\approx n^{\frac{1}{p}}$. It follows by Lemma~\ref{lemmabidemlite} that $\lambda^*\left(n\right)\approx n^{\frac{1}{p'}}$, where $\frac{1}{p}+\frac{1}{p'}=1$. If $p=\infty$, then $\lambda$ is bounded, which entails by unconditionality and seminormalization that $\EB$ is  equivalent to the canonical unit vector basis of $\mathtt{c}_0$. Similarly, if $p=1$, then $p'=\infty$, so $\lambda^*$ is bounded; as $\EB^*$ is also seminormalized and unconditional, it is equivalent to the canonical unit vector basis of $\mathtt{c}_0$, thus $\EB$ is equivalent to the canonical unit vector basis of $\ell_1$. Finally, if $1<p<\infty$,  then $1<p'<\infty$, and $\EB^*$ is bidemocratic by \cite[Proposition 4.4]{DKKT2003}. Hence, $\EB^{**}$ is bidemocratic. As $\EB$ is equivalent to $\EB^{**}$, \cite[Theorem 3.1]{CDKM2023} gives the desired result. \\
\ref{tqg*partially}
By Lemmas~\ref{lemmatranslates} and~\ref{lemmademoc0}, $\EB$ is democratic. Thus, as in the proof of \ref{bothpartial}, $\EB$ is a Schauder basis. It follows by Corollary~\ref{corollaryuncond2} that $\EB$ is unconditional. Now \cite[Theorem 5.1]{CDKM2023} gives the result. \\
\ref{tqg*dem} This is a special case of \ref{tqg*partially*}\ref{polynom}.\\
 \ref{tqg*partially*} 
If \ref{norming2} holds, then $\EB^*$ is seminormalized by Lemma~\ref{lemmanorming}. Thus, in particular, \ref{polynom} holds. \\
If  \ref{polynom}, \ref{lambdasmall2}, or \ref{phinotbig} holds, then $\EB$ is a Schauder basis by Lemma~\ref{lemmaSchauder3}. \\
Thus, by  Corollary~\ref{corollaryuncond2}, $\EB$ is unconditional, and then we can apply~\ref{tqgpartially*}. 
\end{proof}
To finish this note, we give a characterization of the canonical basis of $\ell_p$ or $\mathtt{c}_0$ combining \eqref{p**} with a condition on the translates of its subsequences.
\begin{proposition}\label{propositionsubseqtoellp}(cf. \cite[Theorems 6.1 and 6.2]{CDKM2023} Let $\EB$ be a minimal system for $\XX$. If $\EB$ satisfies \eqref{p**} and every subsequence of $\EB$ is translation invariant, then $\EB$ is equivalent to the canonical unit vector basis of either $\mathtt{c}_0$ or $\ell_p$ for some $1\le p<\infty$. 
\end{proposition}
\begin{proof}
By Lemma~\ref{lemmatranslationtosubsym}, $\EB$ invariant under spreading. Hence, in particular, it is democratic. As $\lambda\approx \varphi$, Lemma~\ref{lemmaSchauder3} entails that $\EB$ is a Schauder basis, and then the result follows from \cite[Theorem 6.1]{CDKM2023}. 
\end{proof}

\section*{Declaration on AI use}

Parts of the text were dictated in Spanish and written in English by OpenAI's ChatGPT. There was no further usage of AI tools.

\end{document}